\documentclass{amsproc}%
\usepackage{amsmath}
\usepackage{graphicx}%
\usepackage{amsfonts}%
\usepackage{amssymb}
\newtheorem{theorem}{Theorem}[section]
\newtheorem{corollary}[theorem]{Corollary}
\newtheorem{lemma}[theorem]{Lemma}
\newtheorem{proposition}[theorem]{Proposition}
\theoremstyle{definition}

\newtheorem{remark}[theorem]{Remark}

\newtheorem{case}{Case}
\theoremstyle{remark}

\newtheorem*{acknowledgements}{Acknowledgements}
\numberwithin{equation}{section}
\newlength{\displayboxwidth}
\newlength{\pairwidth}
\newcounter{enumlink}
\def\openone%{\hbox{\upshape \small1\kern-3.3pt\normalsize1}}
{\mathchoice
{\hbox{\upshape \small1\kern-3.3pt\normalsize1}}
{\hbox{\upshape \small1\kern-3.3pt\normalsize1}}
{\hbox{\upshape \tiny1\kern-2.3pt\SMALL1}}
{\hbox{\upshape \Tiny1\kern-2pt\tiny1}}}
\newbox\ipbox
\newcommand{\ip}[2]{\left\langle #1\mathrel{\mathchoice
{\setbox\ipbox=\hbox{$\displaystyle \left\langle\mathstrut #1#2\right\rangle$}
\vrule height\ht\ipbox width0.25pt depth\dp\ipbox}
{\setbox\ipbox=\hbox{$\textstyle \left\langle\mathstrut #1#2\right\rangle$}
\vrule height\ht\ipbox width0.25pt depth\dp\ipbox}
{\setbox\ipbox=\hbox{$\scriptstyle \left\langle\mathstrut #1#2\right\rangle$}
\vrule height\ht\ipbox width0.25pt depth\dp\ipbox}
{\setbox\ipbox=\hbox{$\scriptscriptstyle \left\langle\mathstrut #1#2\right\rangle$}
\vrule height\ht\ipbox width0.25pt depth\dp\ipbox}
} #2\right\rangle}

\newsavebox{\shaderect}
\renewcommand{\theenumi}{\roman{enumi}}
\renewcommand{\labelenumi}{\textup{(\theenumi )}}
\mathcode`\;="303B
\begin{document}
\title[Wavelet filters and unitary groups]{Wavelet filters and infinite-dimensional unitary groups}
\author{Ola Bratteli}
\address{Department of Mathematics\\
University of Oslo\\
PB 1053 -- Blindern\\
N-0316 Oslo\\
Norway}
\email{bratteli@math.uio.no}
\author{Palle E. T. Jorgensen}
\address{Department of Mathematics\\
The University of Iowa\\
14 MacLean Hall\\
Iowa City, IA 52242-1419\\
U.S.A.}
\email{jorgen@math.uiowa.edu}
\thanks{Research supported by the University of Oslo.}
\begin{abstract}In this paper, we study wavelet filters and their dependence on two numbers,
the scale $N$ and the genus $g$. We show that the wavelet filters, in the
quadrature mirror case, have a harmonic analysis which is based on
representations of the $C^{\ast}$-algebra $\mathcal{O}_{N}$. A main tool in
our analysis is the infinite-dimensional group of all maps $\mathbb{T}%
\rightarrow\mathrm{U}\left(  N\right)  $ (where $\mathrm{U}\left(  N\right)  $
is the group of all unitary $N$-by-$N$ matrices), and we study the extension
problem from low-pass filter to multiresolution filter using this group.
\end{abstract}
\maketitle
\tableofcontents
\listoftables

\section{\label{Int}Introduction}

\setlength{\displayboxwidth}{\textwidth}\addtolength{\displayboxwidth
}{-2\leftmargini}Digital filters are defined by a choice of coefficients which
weight digital time signals from signal processing. It has been known since
the eighties that these coefficients are also those which relate the
operations of scaling and translation on the real line $\mathbb{R}$ in the
construction of multiresolution wavelets; see (\ref{eq1.11}) and
(\ref{eq1.18}). There the digital filters are also known as quadrature mirror
filters and they are fundamental in generating wavelet bases in $L^{2}%
(\mathbb{R})$. They can be determined by various rules which however at times
appear somewhat \emph{ad hoc} in the literature
\cite{JiMi91},
\cite{Dau92},
\cite{BrJo97b}. Here we will describe a
pairwise bijective correspondence between three sets:

\begin{enumerate}
\item \label{Intthings(1)}the wavelet filters,

\item \label{Intthings(2)}a certain infinite-dimensional unitary group (also
called a loop group),

\setcounter{enumlink}{\value{enumi}}
\end{enumerate}

\noindent and finally

\begin{enumerate}
\setcounter{enumi}{\value{enumlink}}

\item \label{Intthings(3)}a specific class of representations on
$L^{2}(\mathbb{T})$ of some relations from $C^{\ast}$-algebra theory. These
last mentioned relations go under the name of the Cuntz relations and are
widely studied in operator theory.
\end{enumerate}

The representations in (\ref{Intthings(3)}) are not
representations of groups, but
rather representations of relations.
These latter relations have been
noted independently in the operator
algebra community and in the
subband filtering community;
see especially \cite{Cun77} and \cite{CMW92a}. In fact,
the present authors stressed in
\cite{BrJo97b} and \cite{BrJo98a} that the digital
filters which go into the analysis
and the synthesis of time
signals with subbands
are special cases of these representations.
In the particular case of quadrature
mirror filters (QMF) there are two bands,
which are called high-pass/low-pass.
In Section \ref{Rep}, we work
out the general case of $N$ bands, and
make explicit the interplay
between the separate viewpoints (\ref{Intthings(1)})--(\ref{Intthings(3)}).
This also has the advantage of
explaining how wavelet packets
(see \cite{CMW92b}) arise very
naturally from the
representation-theoretic method
in (\ref{Intthings(3)}). Band
filtering is described
from an engineering viewpoint in \cite{Vai93}.

The question of reducibility versus irreducibility of the representations in
(\ref{Intthings(3)}) will be central, as will be unitary equivalence of pairs
of representations. Our study of the equivalence question for representations
is motivated by a desire for classifying wavelets and systematically
constructing interesting examples. The question of reducibility of the
representations in turn is motivated by the need for isolating when wavelet
data is minimal in a suitable sense. Wavelets come with some specified
numerical scaling which may be a natural number, $N=2,3,\dots$, or it may be
an expansive integral matrix of some specified dimension $d$. (Here we shall
restrict to $d=1$.) It turns out that the reducibility question for the
representations is closely tied to whether or not the scaling data (in this
case $N$) is minimal in a suitable sense; see Section \ref{Red} below.

We also mention that the algebra $\mathcal{O}_{N}$ is used in the theory of
superselection sectors in quantum field theory \cite{DoRo89}. While
$\mathcal{O}_{N}$ is used there in connection with representation theory, our
viewpoint here is completely different: Doplicher and Roberts use in
\cite{DoRo87} a category of endomorphisms of $\mathcal{O}_{N}$ in deriving a
noncommutative Pontryagin duality theory, also called Tannaka-Krein theory, in
giving an algebraic description of $\hat{G}$ when $G$ is a given compact
(non-abelian) Lie group. The intertwiners for systems of representations of
$G$ in \cite{DoRo89} and \cite{DoRo87} induce the endomorphisms of
$\mathcal{O}_{N}$. Our viewpoint is in a sense the opposite: We identify a
class of representations of $\mathcal{O}_{N}$ which has the structure of a
compact loop group, i.e., the group of all maps from $\mathbb{T}$ into
$\mathrm{U}\left(  N\right)  $ where $\mathrm{U}\left(  N\right)  $ denotes
the group of unitary $N\times N$ matrices.

\section{\label{Bac}Background on wavelet filters and extensions}

To fix notation, let us give a short rundown of the standard multiresolution
wavelet analysis of scale $N$. We follow Section 10 in \cite{BrJo97b}, but see
also \cite{GrMa92}, \cite{Mey87}, \cite{MRV96} and \cite{Hor95}. Define
scaling by $N$ on $L^{2}(\mathbb{R})$ by
\begin{equation}
\left(  U\xi\right)  (x)=N^{-\frac{1}{2}}\xi\left(  N^{-1}x\right)  ,
\label{eq1.1}%
\end{equation}
and translation by $1$ on $L^{2}(\mathbb{R)}$ by
\begin{equation}
\left(  T\xi\right)  (x)=\xi(x-1). \label{eq1.2}%
\end{equation}
A \emph{scaling function} is a function $\varphi\in L^{2}(\mathbb{R})$, such
that if $\mathcal{V}_{0}$ is the closed linear span of all translates
$T^{k}\varphi$, $k\in\mathbb{Z}$, then $\varphi$ has the following four
properties:%
\begin{equation}
\begin{minipage}[t]{\displayboxwidth}\raggedright$\{T^{k}\varphi;k\in
\mathbb{Z}\}$ is an orthonormal set in $L^{2}(\mathbb{R})$; \end{minipage}
\label{eq1.3}%
\end{equation}%
\begin{equation}
\begin{minipage}[t]{\displayboxwidth}\raggedright$U\varphi\in\mathcal{V}%
_{0}$; \end{minipage} \label{eq1.4}%
\end{equation}%
\begin{equation}
\begin{minipage}[t]{\displayboxwidth}\raggedright$\bigwedge_{n\in\mathbb{Z}%
}U^{n}\mathcal{V}_{0}=\{0\}$; \end{minipage} \label{eq1.5}%
\end{equation}%
\begin{equation}
\begin{minipage}[t]{\displayboxwidth}\raggedright$\bigvee_{n\in\mathbb{Z}%
}U^{n}\mathcal{V}_{0}=L^{2}(\mathbb{R}).$ \end{minipage} \label{eq1.6}%
\end{equation}
The simplest example of a scaling function is the Haar function $\varphi$,
which is the characteristic function of the interval $[0,1]$.

By (\ref{eq1.3}) we may define an isometry
\begin{equation}
\mathcal{F}_{\varphi}\colon\mathcal{V}_{0}\longrightarrow L^{2}(\mathbb{T}%
)\colon\xi\longrightarrow m \label{eq1.7}%
\end{equation}
as follows:%
\begin{equation}
\setlength{\arraycolsep}{0pt}%
\begin{array}
[c]{cl}%
\xi(\,\cdot\,) & {}=\sum_{n}b_{n}\varphi(\,\cdot\,-n)\smallskip\\
\llap{$\scriptstyle\mathcal{F}_{\varphi}$}\downarrow\smallskip & \\
m(t) & {}=m(e^{-it})=\sum_{n}b_{n}e^{-int}.
\end{array}
\label{eq1.8}%
\end{equation}
Then
\begin{equation}
\hat{\xi}(t)=m(t)\hat{\varphi}(t) \label{eq1.9}%
\end{equation}
where $\xi\mapsto\hat{\xi}$ is the Fourier transform. If $\xi\in
\mathcal{V}_{-1}=U^{-1}\mathcal{V}_{0}$, then $U\xi\in\mathcal{V}_{0}$, so we
may define%
\begin{equation}
m_{\xi}=\mathcal{F}_{\varphi}\left(  U\xi\right)  \in L^{2}(\mathbb{T}),
\label{eq1.10}%
\end{equation}
and then
\begin{equation}
\sqrt{N}\hat{\xi}(Nt)=m_{\xi}(t)\hat{\varphi}(t)\text{.} \label{eq1.11}%
\end{equation}
In particular, defining
\begin{equation}
m_{0}(t)=m_{\varphi}(t) \label{eq1.12}%
\end{equation}
we note that condition (\ref{eq1.3}) is equivalent to
\begin{equation}
\operatorname{PER}\left(  \left|  \hat{\varphi}\right|  ^{2}\right)
(t):=\sum_{k}\left|  \hat{\varphi}(t+2\pi k)\right|  ^{2}=(2\pi)^{-1}%
\text{,\qquad constancy }\mathrm{a.e.}\;t\in\mathbb{R}\text{,} \label{eq1.13}%
\end{equation}
which in turn implies
\begin{equation}
\sum_{k=0}^{N-1}\left|  m_{0}(t+2\pi k/N)\right|  ^{2}=N. \label{eq1.14}%
\end{equation}
(Note that (\ref{eq1.14}) does \emph{not} conversely imply (\ref{eq1.13}).
Condition (\ref{eq1.14}) merely implies that $\varphi$ defines a tight frame.
In the representation theory in Section \ref{Rep} this distinction does not
play a role. Only the unitarity of $M(z)$ in (\ref{eq1.22}) is relevant. See
\cite[Theorem 3.3.6]{Hor95}. Similarly, the representation theory in Section
\ref{Rep} applies in cases more general than the ones tied to the
$L^{2}(\mathbb{R})$ wavelets via (\ref{eq1.12}) and (\ref{eq1.18}). In the
latter case one also has the low-pass condition $m_{0}(1)=\sqrt{N}$. In the
loop representation (\ref{eq1.23})--(\ref{eq1.24}), this is equivalent to
$A_{0,j}(1)=1/\sqrt{N}$ for $j=0,\ldots,N-1$. These conditions involve the
evaluation of the functions $m_{0}(\,\cdot\,)$ or $A_{0,j}(\,\cdot\,)$ at
$z=1$. But $z=e^{-i\omega}$, where $\omega$ is the frequency variable, so the
conditions are the low-pass conditions for $\omega\sim0.$)

If $\xi,\eta\in U^{-1}\mathcal{V}_{0}$, then we have the equivalence:
\begin{equation}%
\begin{array}
[c]{c}%
\xi\perp T^{k}\eta\text{ for all }k\in\mathbb{Z}\text{\smallskip}\\
\Updownarrow\smallskip\\
\sum_{k=0}^{N-1}\overline{m_{\xi}\left(  t+2\pi k/N\right)  }m_{\eta}\left(
t+2\pi k/N\right)  =0\text{ for almost all }t\in\mathbb{R}\text{.}%
\end{array}
\label{eq1.15}%
\end{equation}
If $\xi\in U^{-1}\mathcal{V}_{0}$, then
\begin{equation}%
\begin{array}
[c]{c}%
\xi(\,\cdot\,-k),\;k\in\mathbb{Z}\text{, is an orthonormal set\smallskip}\\
\Updownarrow\smallskip\\
\sum_{k=0}^{N-1}\left|  m_{\xi}(t+2\pi k/N)\right|  ^{2}=N.
\end{array}
\label{eq1.16}%
\end{equation}
With the \emph{low-pass filter} $m_{0}$ already given, now choose
\emph{high-pass filters} $m_{1},\ldots,m_{N-1}$ in $L^{2}(\mathbb{T})$ such
that
\begin{equation}
\sum_{k=0}^{N-1}\overline{m_{i}\left(  t+2\pi k/N\right)  }m_{j}(t+2\pi
k/N)=\delta_{i,j}N \label{eq1.17}%
\end{equation}
for almost all $t\in\mathbb{R}$, and define \emph{wavelet functions} $\psi
_{1},\ldots,\psi_{N-1}$ by
\begin{equation}
\sqrt{N}\hat{\psi}_{i}(Nt)=m_{i}(t)\hat{\varphi}(t). \label{eq1.18}%
\end{equation}
Using (\ref{eq1.15})--(\ref{eq1.18}) one then shows that the functions
\begin{equation}
\{T^{k}\psi_{i};k\in\mathbb{Z},\;i=1,\ldots,N-1\} \label{eq1.19}%
\end{equation}
form an orthonormal basis for $\mathcal{V}_{-1}\cap\mathcal{V}_{0}\left.
^{\perp}\right.  $, and thus, using (\ref{eq1.5})--(\ref{eq1.6}), the set%
\begin{equation}
\left\{  U^{n}T^{k}\psi_{i};n,k\in\mathbb{Z},\;i=1,\ldots,N-1\right\}
\label{eq1.20}%
\end{equation}
forms an orthonormal basis for $L^{2}(\mathbb{R})$. Concretely, the functions
in (\ref{eq1.20}) are
\begin{equation}
U^{n}T^{k}\psi_{i}(x)=N^{-\frac{n}{2}}\psi_{i}(N^{-n}x-k). \label{eq1.21}%
\end{equation}
Thus, reformulating (\ref{eq1.17}), orthonormality of $\{U^{n}T^{k}\psi_{i}\}$
is equivalent to unitarity of the matrix
\begin{equation}
M(z)=\frac{1}{\sqrt{N}}\left(
\begin{array}
[c]{cccc}%
\vphantom{\vdots}m_{0}(z) & m_{0}(\rho z) & \cdots & m_{0}\left(  \rho
^{N-1}z\right) \\
\vphantom{\vdots}m_{1}(z) & m_{1}(\rho z) & \cdots & m_{1}\left(  \rho
^{N-1}z\right) \\
\vdots &  &  & \\
\vphantom{\vdots}m_{N-1}(z) & m_{N-1}(\rho z) & \cdots & m_{N-1}\left(
\rho^{N-1}z\right)
\end{array}
\right)  \label{eq1.22}%
\end{equation}
for almost all $z\in\mathbb{T}$, where $\rho=e^{\frac{2\mathbb{\pi}i}{N}}$ (so
it is enough to consider $z=e^{ix}$, $0\leq x<\frac{2\mathbb{\pi}}{N}$). Here
$M$ is a function from $\mathbb{T}$ into $\mathrm{U}(N)$ of a special kind. It
will be convenient for our analysis to consider general functions from
$\mathbb{T}$ into $\mathrm{U}(N)$. To this end, we do a Fourier analysis over
the cyclic group $\mathbb{Z}/N\mathbb{Z}$, and we introduce%
\begin{equation}
A_{i,j}(z)=\frac{1}{N}\sum_{w^{N}=z}w^{-j}m_{i}(w), \label{eq1.23}%
\end{equation}
and the inverse transform
\begin{equation}
m_{i}(z)=\sum_{j=0}^{N-1}z^{j}A_{i,j}(z^{N}). \label{eq1.24}%
\end{equation}
These transforms are also the key to the correspondences (\ref{Intthings(1)}%
)~$\Leftrightarrow$ (\ref{Intthings(2)})~$\Leftrightarrow$ (\ref{Intthings(3)}%
) mentioned in the Introduction. Then (\ref{eq1.24}) may be summarized as%
\begin{equation}
M(z)=A(z^{N})\frac{1}{\sqrt{N}}\left(
\begin{array}
[c]{cccc}%
1 & 1 & \cdots & 1\\
z & \rho z & \cdots & \rho^{N-1}z\\
\vdots & \vdots &  & \vdots\\
z^{N-1} & \rho^{N-1}z^{N-1} & \cdots & \rho^{(N-1)^{2}}z^{N-1}%
\end{array}
\right)  \text{;} \label{eq1.25}%
\end{equation}
and the requirement that $M(z)$ is unitary for all $z\in\mathbb{T}$ is now
equivalent to $A(z)$ being unitary for all $z\in\mathbb{T}$. But $A$ is an
arbitrary loop$.$ (A \emph{loop} is by definition a function from $\mathbb{T}$
into $\mathrm{U}\left(  N\right)  $. See the introduction to Section
\ref{Dec}.) Hence the formulas (\ref{eq1.23})--(\ref{eq1.24}) represent the
bijection (\ref{Intthings(1)})--(\ref{Intthings(2)}) alluded to in the
Introduction above.

One main problem considered in this paper is the problem whether the high-pass
filters $m_{1},\ldots,m_{N-1}$ selected such that (\ref{eq1.17}) is valid can
be chosen to be ``nice'' functions when the low-pass filter $m_{0}$ is
``nice'', where ``nice'' may mean continuous, differentiable, polynomial, etc.

Using (\ref{eq1.23}), this problem amounts to choosing a unitary matrix $A(z)$
once the first row of $A(z)$ is given as a unit row vector $a(z)$. Ideally,
one would like to use a selection map $F$ from the unit sphere $S^{2N-1}$ of
$\mathbb{C}^{N}$ into $N$-dimensional orthogonal frames in $\mathbb{C}^{N}$,
i.e.,
\begin{equation}
\mathbf{F}(\mathbf{x})=\left(
\begin{array}
[c]{c}%
\mathbf{F}_{0}(\mathbf{x})\\
\vdots\\
\mathbf{F}_{N-1}(\mathbf{x})
\end{array}
\right)  , \label{eq3.7}%
\end{equation}
where the vectors $\mathbf{F}_{0}(\mathbf{x}),\ldots,\mathbf{F}_{N-1}%
(\mathbf{x})$ form an orthonormal basis for $\mathbb{C}^{N}$ for each
$\mathbf{x}\in S^{2N-1}$ and such that $\mathbf{F}_{0}(\mathbf{x})=\mathbf{x}$
for each $\mathbf{x}$, and define%
\begin{equation}
A(z)=\left(
\begin{array}
[c]{c}%
a(z)\\
\mathbf{F}_{1}(a(z))\\
\vdots\\
\mathbf{F}_{N-1}(a(z))
\end{array}
\right)  . \label{eq3.8}%
\end{equation}
For example, when $N=2$ one may use Daubechies's choice
\begin{equation}
\mathbf{F}_{1}(x_{1},x_{2})=(-\bar{x}_{2},\bar{x}_{1}). \label{eq3.9}%
\end{equation}
In general one can of course find measurable functions with these properties,
but it is remarkable that one cannot choose them to be continuous when $N>2$.
Following
Theorem 4.1 in \cite{JiSh94} and
Remark 10.2 in \cite{BrJo97b}, if $S^{n-1}$ is the unit sphere in
$\mathbb{R}^{n}$, then a theorem of Adams \cite{Ada62} says that the highest
number of pointwise linearly independent vector fields that can be defined on
$S^{n-1}$ is $\rho(n)-1$, where the function $\rho(n)$ is defined as follows:
Let $b$ be the multiplicity of $2$ in the prime decomposition of $n$; write
$b=c+4d$ where $c\in\{0,1,2,3\}$ and $d\in\{0,1,2,\ldots\}$; and put
$\rho(n)=2^{c}+8d$. Thus one verifies $\rho(2)=2$, $\rho(4)=4$, $\rho(8)=8$,
$\rho(16)=9$, and in general $\rho(n)<n$ if $n\notin\{2,4,8\}$. Hence the only
possibilities of finding a continuous selection function $F$ are when $N=2$ or
$N=4$. When $N=2$ one can use Daubechies's selection function (\ref{eq3.9}).
When $N=4$ it is tempting to use the Cayley numbers to construct a selection
function, but we will see that this is impossible. Note that if
$\mathrm{U}(N)$ is replaced by the real orthogonal group $\mathrm{O}(N)$, a
selection function can be found if and only if $N\in\{2,4,8\}$ by using the
complex numbers, the quaternions and the Cayley numbers, respectively. For
example, for $N=4$, the matrix
\begin{equation}
\left(
\begin{array}
[c]{rrrr}%
z_{1} & z_{2} & z_{3} & z_{4}\\
-z_{2} & z_{1} & -z_{4} & z_{3}\\
-z_{3} & z_{4} & z_{1} & -z_{2}\\
-z_{4} & -z_{3} & z_{2} & z_{1}%
\end{array}
\right)  \label{eq3.10}%
\end{equation}
is in $\mathrm{O}(4)$ whenever the first row is a unit vector in
$\mathbb{R}^{4}$. Using the identification of $\mathbb{C}$ as $2\times2$ real
matrices%
\begin{equation}
z_{1}+iz_{2}\approx\left(
\begin{array}
[c]{rr}%
z_{1} & z_{2}\\
-z_{2} & z_{1}%
\end{array}
\right)  , \label{eq3.11}%
\end{equation}
this corresponds to Daubechies's selection function (\ref{eq3.9}).
If $N=4$, a selection function does
not exist in the complex case, by Proposition 3.6 in
\cite{Jam76}. Specializing
to the case $n=k=N$ in this proposition, it
says that, if the space of orthonormal
$N$-frames over $\mathbb{C}^{N}$ admits a cross-section
over the unit sphere in $\mathbb{C}^{N}$, and if $p$
is any prime, then $N$ is divisible by the
smallest power of $p$ exceeding $\left( N-1\right) /\left( p-1\right) $.
If $N=4$, this condition fails for $p=3$
(while it is fulfilled for all $p$ if $N=2$).
(Note that a direct application of Adams's
theorem only gives nonexistence of the
cross-section for $N>4$, as is also
indicated in the proof of Theorem 4.1 in
\cite{JiSh94}.)
(Note also that the discussion above settles the question
implied by footnote 3 to the remark before Lemma 3.2.3
in \cite{Hor95}.)
The claim made in \cite[Remark 10.2]{BrJo97b} that there is such a
selection
function for $N=8$ is definitely erroneous. What is true is that the matrix%
\begin{equation}
\left(
\begin{array}
[c]{rrrrrrrr}%
z_{1} & z_{2} & z_{3} & z_{4} & z_{5} & z_{6} & z_{7} & z_{8}\\
-z_{2} & z_{1} & -z_{4} & z_{3} & -z_{6} & z_{5} & z_{8} & -z_{7}\\
-z_{3} & z_{4} & z_{1} & -z_{2} & -z_{7} & -z_{8} & z_{5} & z_{6}\\
-z_{4} & -z_{3} & z_{2} & z_{1} & -z_{8} & z_{7} & -z_{6} & z_{5}\\
-z_{5} & z_{6} & z_{7} & z_{8} & z_{1} & -z_{2} & -z_{3} & -z_{4}\\
-z_{6} & -z_{5} & z_{8} & -z_{7} & z_{2} & z_{1} & z_{4} & -z_{3}\\
-z_{7} & -z_{8} & -z_{5} & z_{6} & z_{3} & -z_{4} & z_{1} & z_{2}\\
-z_{8} & z_{7} & -z_{6} & -z_{5} & z_{4} & z_{3} & -z_{2} & z_{1}%
\end{array}
\right)  \label{eq3.11bis}%
\end{equation}
is in $\mathrm{O}\left(  8\right)  $ whenever the first row is a unit vector
in $\mathbb{R}^{8}$. However, trying to convert this to an element of
$\mathrm{U}\left(  4\right)  $ by the same trick as above leads to%
\begin{equation}
\left(
\begin{array}
[c]{rrrr}%
z_{1} & z_{2} & z_{3} & z_{4}\\
-\overline{z_{2}} & \overline{z_{1}} & -z_{4} & z_{3}\\
-\overline{z_{3}} & z_{4} & \overline{z_{1}} & -z_{2}\\
-z_{4} & -\overline{z_{3}} & \overline{z_{2}} & z_{1}%
\end{array}
\right)  , \label{eq3.11ter}%
\end{equation}
but this is \emph{not} in general unitary when the first row is a unit vector
in $\mathbb{C}^{4}$.
For an interesting short account of the matrices
(\ref{eq3.10})--(\ref{eq3.11bis}), see
\cite{Tau71}. The similar selection problem in the context of splines
and higer-dimensional spaces $\mathbb{R}^{s}$ has been
treated in \cite{RiSh91} and \cite{RiSh92}; see in particular
pp.~142--145 in \cite{RiSh91}.

However, there is a completely different method of selecting $A(z)$ using more
global properties of the first row $a(z)$ which applies in the case that $a$
is a polynomial in $z$. This special case is interesting because it
corresponds to compactly supported scaling functions $\varphi$: the scaling
function $\varphi$ constructed from $m_{0}$ by cascade approximation has
support in $[0,Ng-1]$ if and only if the low-pass filter $m_{0}(z)$ is a
polynomial of degree at most $g-1$, and satisfies $m_{0}(1)=\sqrt{N}$. See
\cite{BrJo97b}, \cite{BrJo98a}, \cite{BrJo99b}, \cite{BEJ00}, \cite{ReWe98},
and \cite{Dau92}. In this context there is a method of extending $a(z)$ to
$A(z)$ which is independent of the selection theory, and $A(z)$ may even be
taken to be a polynomial in $z$ (with coefficients in $\mathcal{B}%
(\mathbb{C}^{N})$) of degree equal to the degree of $a$. The embryonic form of
this method is due to Pollen \cite{Pol90}, \cite{Pol92}, but has been
developed further: see, e.g., \cite{ReWe98}, \cite{Vai93}. We will give an
account of the method in operator-theoretic form, and give the results in
Proposition \ref{proposition2.2} and Theorem \ref{ThmExt.1}. Combining Theorem
\ref{ThmExt.1} with (\ref{eq1.25}) we therefore obtain

\begin{theorem}
\label{theorem1.1}Let $N\in\{2,3,\ldots\}$ and let $\varphi$ be a scaling
function with support in $[0,Ng-1]$ where $g\in\{1,2,3,\ldots\}$. Then one can
use multiresolution wavelet analysis to find associated wavelet functions
$\psi_{1},\ldots,\psi_{N-1}$ also having support in $[0,Ng-1]$.
\end{theorem}

\begin{proof}
The conditions on $\varphi$, $\psi_{1},\;\ldots$, $\psi_{N-1}$ correspond to
the functions $m_{0}$, $m_{1},\;\ldots$, $m_{N-1}$ being polynomials of degree
$Ng-1$, \cite{Hor95}, hence Theorem \ref{theorem1.1} follows from Theorem
\ref{ThmExt.1} and (\ref{eq1.25}). More details are given in Corollary
\ref{CorExt.2} and its proof.
\end{proof}

In Sections \ref{Rep} and \ref{Red} we will consider the representations of
the Cuntz algebra $\mathcal{O}_{N}$ defined by $m_{0},m_{1},\ldots,m_{N-1}$
when these are polynomials, and thus extend results from \cite{BEJ00}.

\section{\label{Dec}Decomposition of polynomial loops into linear factors}

We define a \emph{loop} to be a continuous function from the circle
$\mathbb{T}$ into the compact group $\mathrm{U}(N)$ of unitary $N\times N$
matrices, where $N=1,2,\dots$. The set $C(\mathbb{T},\mathrm{U}(N))$ of loops
has a natural group structure under pointwise multiplication, and this group
is called the \emph{loop group} \cite{PrSe86}. We say that a loop is
\emph{polynomial} if its matrix elements are polynomials in the basic variable
$z\in\mathbb{T}=\{w\in\mathbb{C};\left|  w\right|  =1\}$. The set
$\mathcal{P}(\mathbb{T},\mathrm{U}(N))$ of polynomial loops then forms a
semigroup in $C(\mathbb{T},\mathrm{U}(N))$ called the \emph{polynomial loop}
\emph{semigroup}. (Note that the formulas (\ref{eq1.23})--(\ref{eq1.24}) for
the bijection between the set of wavelet filters and the loop group are valid
also in the wider category of $L^{\infty}$-functions, i.e., when each
$m_{i}^{(A)}\in L^{\infty}(\mathbb{T})$ and each $A_{i,j}\in L^{\infty
}(\mathbb{T})$. This is the generality of \cite{BrJo97b}.) Note that if $U$ is
a given loop, then $U$ has a class in the $K_{1}$-group of $C(\mathbb{T}%
,M_{N})$, i.e., $K_{1}\left(  C\left(  \mathbb{T},M_{N}\right)  \right)
\cong\mathbb{Z}$, see \cite{Bla86}. This class is called the McMillan degree
in the wavelet literature \cite{ReWe98}, and we will denote it by
$K_{1}\left(  U\right)  $. It can be computed as the winding number of the
map
\begin{equation}
\mathbb{T}\ni z\longrightarrow\det\left(  U(z)\right)  \in\mathbb{T}\text{.}
\label{eq2.1}%
\end{equation}
If $U$ is in the loop semigroup, the integer $K_{1}(U)$ is necessarily
nonnegative.
This follows from Lemma \ref{lemma2.1}, below,
applied to $u\left( z\right) =\det\left( U\left( z\right) \right) $.
The map $U\rightarrow K_{1}(U)$ is a group homomorphism from
$C(\mathbb{T},\mathrm{U}(N))$ onto $\mathbb{Z}$. We need an elementary lemma
\cite{ReWe98}.

\begin{lemma}
\label{lemma2.1}If $u\colon\mathbb{T}\rightarrow\mathbb{T}$ is a polynomial,
then $u$ is a monomial.
\end{lemma}

\begin{proof}
We have $u=\sum_{k=0}^{n}c_{k}z^{k}$ for suitable coefficients $c_{k}%
\in\mathbb{Z}$. By multiplying $u$ by a nonpositive power of $z$ we may assume
$c_{0}\neq0$. But then $1=u(z)\overline{u(z)}=c_{n}\overline{c_{0}}z^{n}+$
powers of $z$ between $-(n-1)$ and $(n-1)$ $+\overline{c_{n}}c_{0}z^{-n}$ and
hence $c_{n}=0$. By induction, $u(z)=c_{0}$.
\end{proof}

We now use this to show that any polynomial loop $U$ has a decomposition into
linear factors \cite[pp.~60--61]{ReWe98}.
The idea in the proof of Lemma \ref{lemma2.3} goes at least back to
\cite{Pol90}
and \cite[Lemma~3.3]{LLS96}.

\begin{proposition}
\label{proposition2.2}Let $U$ be a polynomial loop in $\mathcal{P}%
(\mathbb{T},\mathrm{U}(N))$. The following conditions are equivalent.%
\begin{equation}
\begin{minipage}[t]{\displayboxwidth}\raggedright$K_{1}%
(U)=d$. \end{minipage}\label{eq2.2}
\end{equation}
\begin{equation}
\begin{minipage}[t]{\displayboxwidth}\raggedright
There exist one-dimensional projections $P_{1},\dots,P_{d}$ on $\mathbb{C}%
^{N}$ and a $V\in\mathrm{U}(N)$ such that $$ U(z)=U_{1}(z)U_{2}(z)\cdots
U_{d}(z)V,$$ where $$ U_{i}(z)=zP_{i}+(1-P_{i}) $$ for $i=1,\dots
,d$ and $z\in\mathbb{T}$. \textup{(}If $d=0$, then $U(z)=V$.\textup{)}
\end{minipage}
\label{eq2.3}
\end{equation}
\end{proposition}

\begin{proof}
Since $K_{1}(U_{i})=1$ and $K_{1}(V)=0$, the implication (\ref{eq2.3})
$\Rightarrow$ (\ref{eq2.2}) is obvious. We prove the other implication by
using induction with respect to $d$. To this end, we will use the following
reduction lemma.\renewcommand{\qed}{}

\begin{lemma}
\label{lemma2.3}Assume that
\begin{equation}
U(z)=A_{0}+A_{1}z+\cdots+A_{k}z^{k} \label{eq2.4}%
\end{equation}
defines a $U\in\mathcal{P}(\mathbb{T},\mathrm{U}(N))$, where $A_{k}\neq0$ and
$k\geq1$. Then $U$ has a unique decomposition%
\begin{equation}
U(z)=((1-Q)+zQ)(B_{0}+B_{1}z+\cdots+B_{k-1}z^{k-1}), \label{eq2.5}%
\end{equation}
where $Q$ is the projection onto the range $A_{k}\mathbb{C}^{N}$ of $A_{k}$
\textup{(}and then $(1-Q)+zQ$ and $B_{0}+B_{1}z+\cdots+B_{k-1}z^{k-1}$ are in
$\mathcal{P}(\mathbb{T},\mathrm{U}(N))$\textup{).}
\end{lemma}
\end{proof}

\begin{proof}
Let $Q$ be the projection onto $A_{k}\mathbb{C}^{N}$. We have
\begin{align}
1  &  =U(z)^{\ast}U(z)\label{eq2.6}\\
&  =(A_{0}^{\ast}+A_{1}^{\ast}z^{-1}+\cdots+A_{k}^{\ast}z^{-k})(A_{0}%
+A_{1}z+\cdots+A_{k}z^{k})\nonumber\\
&  =A_{k}^{\ast}A_{0}z^{-k}+\text{ terms in higher powers of }z\text{.}%
\nonumber
\end{align}
Hence%
\begin{equation}
A_{k}^{\ast}A_{0}=0, \label{eq2.7}%
\end{equation}
and thus, since $Q$ is the projection onto the range of $A_{k}$, we obtain
\begin{equation}
QA_{0}=0\text{\qquad and\qquad}(1-Q)A_{k}=0\text{.} \label{eq2.8}%
\end{equation}
Now put
\begin{align}
W(z)  &  =((1-Q)+zQ)^{-1}U(z)\label{eq2.9}\\
&  =((1-Q)+z^{-1}Q)(A_{0}+\cdots+z^{k}A_{k})\nonumber\\
&  =z^{-1}QA_{0}\nonumber\\
&  \qquad+((1-Q)A_{0}+QA_{1})\nonumber\\
&  \qquad+\cdots\nonumber\\
&  \qquad+z^{k}(1-Q)A_{k}\text{\quad(ordered by increasing powers).}\nonumber
\end{align}
It follows from (\ref{eq2.8}) that $W(z)$ has the form
\begin{equation}
W(z)=B_{0}+B_{1}z+\cdots+B_{k-1}z^{k-1}, \label{eq2.10}%
\end{equation}
and (\ref{eq2.5}) follows. (In fact, $B_{0}=\left(  1-Q\right)  A_{0}%
+QA_{1}=A_{0}+QA_{1}$, $\dots$, and $B_{k-1}=\left(  1-Q\right)  A_{k-1}%
+A_{k}$.) By our choice of $Q$, the decomposition (\ref{eq2.5}) is unique by
(\ref{eq2.8}).
\end{proof}

\begin{proof}
[End of the proof of Proposition \textup{\ref{proposition2.2}}]We will prove
the proposition by induction with respect to the McMillan index $d$. Assume
first that $d\geq1$ and that the proposition holds for all indices $\leq d-1$,
and suppose that
\begin{equation}
U(z)=A_{0}+A_{1}z+\cdots+A_{k}z^{k} \label{eq2.11}%
\end{equation}
has index $d$. But by Lemma \ref{lemma2.3}, $U$ has a decomposition
\begin{equation}
U(z)=((1-Q)+zQ)W(z), \label{eq2.12}%
\end{equation}
where $Q$ is a nonzero projection. But then
\begin{align}
d  &  =K_{1}(U)=K_{1}((1-Q)+zQ)+K_{1}(W(z))\label{eq2.13}\\
&  =\dim(Q)+K_{1}(W(z)),\nonumber
\end{align}
so $K_{1}(W(z))\leq d-1$.
(We define $\dim Q=\mathop{\mathrm{rank}} Q$ when $Q$ is a projection.)
We may then apply the induction hypothesis to
$W(z)$. Finally, there exist $\dim Q$ one-dimensional projections
$P_{1},\ldots,P_{\dim(Q)}$ such that $Q=P_{1}+P_{2}+\cdots+P_{\dim(Q)}$, and
then
\begin{equation}
((1-Q)+zQ)=\prod_{n=1}^{\dim Q}((1-P_{n})+zP_{n}). \label{eq2.14}%
\end{equation}
It follows that $U$ has the form (\ref{eq2.3}).

Finally, if $K_{1}(U)=0$ and $U$ has the form (\ref{eq2.11}), there are two
possibilities:\renewcommand{\theenumi}{\arabic{enumi}} \renewcommand
{\labelenumi}{\textup{\theenumi.}}

\begin{enumerate}
\item \label{proposition2.2proof(1)}$k=0$ and $U$ is a constant unitary
matrix. Then $U$ already has the form in (\ref{eq2.3}) with $d=0$.

\item \label{proposition2.2proof(2)}$k>0$. Then one may apply Lemma
\ref{lemma2.3} to write $U(z)=((1-Q)+zQ)W(z)$ with $Q$ a nonzero projection,
but as%
\[
0=K_{1}(U)=\dim Q+K_{1}(W)
\]
and $K_{1}(W)\geq0$, this is impossible.
\end{enumerate}

\noindent
This ends the proof of Proposition \ref{proposition2.2}.
The last argument in the induction chain
is simpler if we use induction with
respect to $k$ rather than induction with
respect to $d$.
\end{proof}

\section{\label{Ext}Extensions of low-pass polynomial filters to high-pass filters}

In this section we will prove Theorem \ref{theorem1.1}. This theorem follows
from Corollary \ref{CorExt.2}, below. Theorem \ref{ThmExt.1} states that every
vector $\alpha$ in $\mathbb{C}^{Ng}$ which satisfies a certain orthogonality
condition is the first row of coefficients of some element of the loop group.

\begin{theorem}
\label{ThmExt.1}Let $\alpha=\left(  \alpha_{0},\alpha_{1},\dots,\alpha
_{g-1}\right)  $ be $g$ row vectors in $\mathbb{C}^{N}$. The following
conditions are equivalent.%
\begin{equation}
\begin{minipage}[t]{\displayboxwidth}\raggedright
The vectors satisfy the relations $$ \sum_{i}\ip{\alpha_{i}}{\alpha_{i+j}%
}=0 \text{\qquad for }j\neq0, $$ and $$ \sum_{i}\ip{\alpha_{i}}{\alpha_{i}%
}=1, $$ where we use the convention that $\alpha_{i}=0$ if $i<0$ or $i\geq
g$. \end{minipage}\label{eqExt.1}
\end{equation}
\begin{equation}
\begin{minipage}[t]{\displayboxwidth}\raggedright
There exists a polynomial loop $$ A\left( z\right) \in\mathcal{P}%
\left( \mathbb{T},\mathrm{U}\left( N\right) \right
) $$ of degree $g-1$ such that the first row of $A\left( z\right
) $ is $$ \sum_{i=0}^{g-1}z^{i}\alpha_{i}. $$ \end{minipage}\label{eqExt.2}
\end{equation}
\end{theorem}

\begin{proof}
The implication (\ref{eqExt.2}) $\Rightarrow$ (\ref{eqExt.1}) follows from
considering the
$\left( 0,0\right) $ matrix entry
of%
\begin{equation}
A\left(  z\right)  A\left(  z\right)  ^{\ast}=1 \label{eqExt.3}%
\end{equation}
in $M_{N}$ for all $z$, i.e., $A\left(  z\right)  A\left(  z\right)  ^{\ast}$
is the constant Laurent polynomial $1$.

The other implication is proved by a very similar method as in the proof of
Proposition \ref{proposition2.2}. Again we use induction with respect to $g$.

If $g=1$, condition (\ref{eqExt.1}) just says that $\alpha_{0}$ is a unit row
vector, and we can find a constant function $A\left(  z\right)  =A$ just by
Gram-Schmidt orthogonalization.

Assume that $g>1$, and that the result has been proved for all smaller $g$. We
may assume that $\alpha_{g-1}\neq0$ (otherwise we are already through by the
induction hypothesis). Let $P$ be the one-dimensional projection onto
$\alpha_{g-1}$, i.e.,%
\begin{equation}
\alpha_{g-1}P=\alpha_{g-1}. \label{eqExt.4}%
\end{equation}
Now define%
\begin{align}
\beta\left(  z\right)   &  =\alpha\left(  z\right)  \left(  1-P+z^{-1}P\right)
\label{eqExt.5}\\
&  =\left(  \alpha_{0}+\alpha_{1}z+\dots+\alpha_{g-1}z^{g-1}\right)  \left(
1-P+z^{-1}P\right) \nonumber\\
&  =z^{-1}\alpha_{0}P\nonumber\\
&  \qquad+\alpha_{0}\left(  1-P\right)  +\alpha_{1}P\nonumber\\
&  \qquad+z\left(  \cdots\right) \nonumber\\
&  \qquad+\dots\nonumber\\
&  \qquad+z^{g-2}\left(  \alpha_{g-1}P+\alpha_{g-2}\left(  1-P\right)  \right)
\nonumber\\
&  \qquad+z^{g-1}\alpha_{g-1}\left(  1-P\right)  .\nonumber
\end{align}
But since $\ip{\alpha_{0}}{\alpha_{g-1}}=0$ by (\ref{eqExt.1}), the $z^{-1}$
term disappears, and (\ref{eqExt.4}) implies that the $z^{g-1}$ term
disappears. Therefore $\beta\left(  z\right)  $ is a polynomial in $z$ of
degree $g-2$. One now verifies from the unitarity of $\left(  1-P+z^{-1}%
P\right)  $ that the coefficient vectors of $\beta$ satisfy the same relations
(\ref{eqExt.1}) as those of $\alpha$. Hence, applying the induction
hypothesis, there exists a polynomial loop $B\left(  z\right)  $ of degree
$g-2$ such that the first row of $B\left(  z\right)  $ is $\beta\left(
z\right)  $. Thus, putting%
\begin{equation}
A\left(  z\right)  =B\left(  z\right)  \left(  1-P+zP\right)  ,
\label{eqExt.6}%
\end{equation}
it follows from (\ref{eqExt.5}) that the first row of $A\left(  z\right)  $ is
$\alpha\left(  z\right)  $. This completes the induction, and thus the proof
of Theorem \ref{ThmExt.1}.
\end{proof}

Let us restate Theorem \ref{theorem1.1} and its proof in terms of the low-pass
filter $m_{0}(z)$:

\begin{corollary}
\label{CorExt.2}Let $m_{0}$ be a polynomial, and let $N\in\left\{
2,3,\dots\right\}  $. Suppose%
\begin{equation}
\left|  m_{0}\left(  z\right)  \right|  ^{2}+\left|  m_{0}\left(
z\rho\right)  \right|  ^{2}+\dots+\left|  m_{0}\left(  z\rho^{N-1}\right)
\right|  ^{2}=N \label{eqCorExt.2}%
\end{equation}
for all $z\in\mathbb{T}$, where $\rho=e^{\frac{2\pi i}{N}}$. Then there are
polynomials $\left\{  m_{i};i=1,\dots,N-1\right\}  $
such that the combined system $\left\{  m_{i};0\leq i<N\right\}  $
satisfies the unitarity property of \textup{(\ref{eq1.22}),} or equivalently
\textup{(\ref{eq1.17})} with the convention $z\leftrightarrow e^{it}$. In
other words, every $m_{0}$ may be completed to a quadrature mirror filter system.

Furthermore, when $g\in\mathbb{N}$ is chosen such that
the degree of the polynomial $m_{0}$ is at most
$Ng-1$, then the polynomials $m_{1},\dots ,m_{N-1}$ can
be taken to have degree at most $Ng-1$.
\end{corollary}

\begin{proof}
With $m_{0}$ and $N$ given as stated in the corollary, set%
\begin{equation}
A_{0,j}\left(  z\right)  :=\frac{1}{N}\sum_{w^{N}=z}w^{-j}m_{0}\left(
w\right)  ,\qquad z\in\mathbb{T}. \label{eqCorExtProof.1}%
\end{equation}
Then it follows from (\ref{eqCorExt.2}) that the coefficients $\alpha_{j}$, as
a set of row vectors in $\mathbb{C}^{N}$, satisfy (\ref{eqExt.1}) in Theorem
\ref{ThmExt.1}. Picking then $A\left(  z\right)  =\left(  A_{i,j}\left(
z\right)  \right)  \in\mathcal{P}\left(  \mathbb{T},\mathrm{U}\left(
N\right)  \right)  $ as in (\ref{eqExt.2}), using the theorem, and setting%
\begin{equation}
m_{i}\left(  z\right)  =\sum_{j=0}^{N-1}z^{j}A_{i,j}\left(  z^{N}\right)  ,
\label{eqCorExtProof.2}%
\end{equation}
see (\ref{eq1.24}), we note that these functions satisfy the conclusion in the corollary.

Finally, we note that, as in Remark \ref{RemRepNew.6},
if $m_{0}$ has degree at most $Ng-1$, it follows
from (\ref{eqCorExtProof.1}) that
$A_{0,j}\left( z\right) $ has degree at most $g-1$;
thus $A_{i,j}\left( z\right) $ can be taken to have degree at most $g-1$
by Theorem \ref{ThmExt.1}, and then all
$m_{i}\left( z\right) $ have degree at most
$Ng-1$ by Remark \ref{RemRepNew.6} or (\ref{eqCorExtProof.2}).
\end{proof}

\section{\label{Rep}Representations associated with polynomial filters}

If $m_{0},m_{1},\dots,m_{N-1}$ are complex functions on $\mathbb{T}$ such that
the matrix $M\left(  z\right)  $ in (\ref{eq1.22}) is unitary, one may define
isometries $S_{0},\dots,S_{N-1}$ on $L^{2}\left(  \mathbb{T}\right)  $ by%
\begin{equation}
\left(  S_{i}\xi\right)  \left(  z\right)  =m_{i}\left(  z\right)  \xi\left(
z^{N}\right)  \label{eqRep.1}%
\end{equation}
and then%
\begin{equation}
\left(  S_{i}^{\ast}\xi\right)  \left(  z\right)  =\frac{1}{N}\sum
_{\substack{w\\w^{N}=z}}\bar{m}_{i}\left(  w\right)  \xi\left(  w\right)  ,
\label{eqRep.2}%
\end{equation}
see \cite[(1.16)--(1.17)]{BrJo97b}
and, for $N=2$, \cite[(1.1)]{CMW92a}.
The isometries $S_{0},\dots,S_{N-1}$ then
have mutually orthogonal ranges,%
\begin{equation}
S_{i}^{\ast}S_{j}^{{}}=\delta_{i,j}\openone, \label{eqRep.3}%
\end{equation}
and these ranges add up to $\openone$,%
\begin{equation}
\sum_{i=0}^{N-1}S_{i}^{{}}S_{i}^{\ast}=\openone. \label{eqRep.4}%
\end{equation}
The relations (\ref{eqRep.3})--(\ref{eqRep.4}) are called the \emph{Cuntz
relations} (of order $N$), and the universal $C^{\ast}$-algebra generated by
operators $s_{0},\dots,s_{N-1}$ satisfying these relations is called the
\emph{Cuntz algebra} (of order $N$). It is a simple antiliminal $C^{\ast}%
$-algebra, which means that its space of Hilbert-space realizations cannot be
equipped with a standard Borel structure (see \cite{Cun77}).
In \cite{BrJo97b} and Section \ref{Bac} above,
we identify intertwining operators $W=\mathcal{F}_{\varphi}^{\ast}$
from $L^{2}\left( \mathbb{T}\right) $ onto resolution subspaces in
$L^{2}\left( \mathbb{R}\right) $ which intertwine the isometry $S_{0}$ with
the scaling operator
$f\mapsto N^{-1/2}f\left( x/N\right) $ on the appropriate resolution subspace
in $L^{2}\left( \mathbb{R}\right) $. The family of subspaces
\begin{equation}
S_{i_{1}}S_{i_{2}}\cdots S_{i_{k}}L^{2}\left( \mathbb{T}\right) ,\qquad 
k=0,1,\dots ,\; i_{j}\in \left\{ 0,1,\dots ,N-1\right\} , \label{eqRep.4bis}
\end{equation}
corresponds under $W$ to wavelet packets (see \cite{CMW92b}) in
$L^{2}\left( \mathbb{R}\right) $.
Relations (\ref{eqRep.3})--(\ref{eqRep.4}) make it immediately clear that
the respective projections onto the subspaces (\ref{eqRep.4bis})
are
$S^{}_{i_{1}}\cdots S^{}_{i_{k}}S^{\ast}_{i_{k}}\cdots S^{\ast}_{i_{1}}$.
In particular,
these projections are clearly mutually orthogonal when $k$ is fixed and
multi-indices are varied.

Let us recall some results from \cite{BJKW00} and \cite{BEJ00}. Let
$\mathcal{H}$ be a Hilbert space, and $s_{i}\rightarrow\pi\left(
s_{i}\right)  =S_{i}\in\mathcal{B}\left(  \mathcal{H}\right)  $ a
representation of the Cuntz relations on $\mathcal{H}$. Assume that there
exists a finite-dimensional subspace $\mathcal{K}\subset\mathcal{H}$ with the
properties%
\begin{equation}
S_{i}^{\ast}\mathcal{K}\subset\mathcal{K},\qquad i=0,1,\dots,N-1,
\label{eqRep.5}%
\end{equation}
and%
\begin{equation}
\mathcal{K}\text{ is cyclic, i.e., }\left[  \pi\left(  \mathcal{O}_{N}\right)
\mathcal{K}\right]  =\mathcal{H} \label{eqRep.6}%
\end{equation}
(here $\left[  \pi\left(  \mathcal{O}_{N}\right)  \mathcal{K}\right]  $
denotes the closure of the linear span of all polynomials in\linebreak
$S_{0},\dots,S_{N-1}$ applied to vectors in $\mathcal{K}$). Define operators
$V_{i}\in\mathcal{B}\left(  \mathcal{K}\right)  $ by%
\begin{equation}
V_{i}^{\ast}=S_{i}^{\ast}|_{\mathcal{K}}. \label{eqRep.7}%
\end{equation}
Then%
\begin{equation}
\sum_{i=0}^{N-1}V_{i}^{{}}V_{i}^{\ast}=\openone_{\mathcal{K}}. \label{eqRep.8}%
\end{equation}
Define a map $\sigma$ on $\mathcal{B}\left(  \mathcal{K}\right)  $ by%
\begin{equation}
\sigma\left(  A\right)  =\sum_{i=0}^{N-1}V_{i}^{{}}AV_{i}^{\ast}.
\label{eqRep.9}%
\end{equation}
Then $\sigma$ is a unital completely positive map, and

\begin{theorem}
\label{ThmRep.1}\textup{(\cite{BJKW00})} There is a positive norm-preserving
linear isomorphism between the commutant algebra%
\begin{equation}
\pi\left(  \mathcal{O}_{N}\right)  ^{\prime}=\left\{  A\in\mathcal{B}\left(
\mathcal{H}\right)  ;A\pi\left(  x\right)  =\pi\left(  x\right)  A\text{ for
all }x\in\mathcal{O}_{N}\right\}  \label{eqRep.10}%
\end{equation}
and the fixed-point set%
\begin{equation}
\mathcal{B}\left(  \mathcal{K}\right)  ^{\sigma}=\left\{  A\in\mathcal{B}%
\left(  \mathcal{K}\right)  ;\sigma\left(  A\right)  =A\right\}
\label{eqRep.11}%
\end{equation}
given by%
\begin{equation}
\pi\left(  \mathcal{O}_{N}\right)  ^{\prime}\ni A\longrightarrow PAP,
\label{eqRep.12}%
\end{equation}
where $P$ is the projection of $\mathcal{H}$ onto $\mathcal{K}$. In
particular, $\pi$ is irreducible if and only if $\sigma$ is ergodic.

More generally, if $\mathcal{K}_{1}$, $\mathcal{K}_{2}$ \textup{(}with
corresponding projections $P^{\left(  1\right)  }$ and $P^{\left(  2\right)
}$\textup{)} are $S^{\ast}$-invariant cyclic subspaces for two representations
$\pi_{1}$, $\pi_{2}$ of $\mathcal{O}_{N}$ on $\mathcal{H}_{1}$, $\mathcal{H}%
_{2}$, and%
\begin{equation}
V_{i}^{\left(  j\right)  }=P_{{}}^{\left(  j\right)  }\pi_{j}^{{}}\left(
s_{i}\right)  |_{\mathcal{K}_{j}} \label{eqRep.13}%
\end{equation}
for $j=1,2$, $i=0,\dots,N-1$, define $\rho$ on $\mathcal{B}\left(
\mathcal{K}_{1},\mathcal{K}_{2}\right)  $ by%
\begin{equation}
\rho\left(  A\right)  =\sum_{i}V_{i}^{\left(  2\right)  }AV_{i}^{\left(
1\right)  \,\ast}. \label{eqRep.14}%
\end{equation}
Then there is an isometric linear isomorphism between the set of intertwiners%
\begin{equation}
\left\{  A\in\mathcal{B}\left(  \mathcal{H}_{1},\mathcal{H}_{2}\right)
;A\pi_{1}\left(  x\right)  =\pi_{2}\left(  x\right)  A\text{ for all }%
x\in\mathcal{O}_{N}\right\}  \label{eqRep.15}%
\end{equation}
and the fixed-point set%
\begin{equation}
\left\{  B\in\mathcal{B}\left(  \mathcal{K}_{1},\mathcal{K}_{2}\right)
;\rho\left(  B\right)  =B\right\}  \label{eqRep.16}%
\end{equation}
given by%
\begin{equation}
A\longrightarrow B=P^{\left(  2\right)  }AP^{\left(  1\right)  }.
\label{eqRep.17}%
\end{equation}
\end{theorem}

We argued in \cite[(4.14)--(4.18)]{BEJ00}
that if $m_{0},\dots,m_{N-1}$ are polynomials in the
circle variable $z$, $N=2$,
and $\mathcal{H}=L^{2}\left(  \mathbb{T}\right)  $, then
such a finite-dimensional subspace $\mathcal{K}$ exists, having dimension $Ng$
and spanned by $1,z^{-1},z^{-2},\dots,z^{-Ng+1}$. In this section we will
extend this result to $N>2$,
and we will see in Proposition \ref{ProRep.5} that
we can do slightly better than what the
$N=2$ result indicates.
Let us
assume from now on that the loop group element $A\left(  z\right)  $ in
(\ref{eq1.23})--(\ref{eq1.25}) is a polynomial of degree $g-1$, so that
$m_{0}\left(  z\right)  ,\dots,m_{N-1}\left(  z\right)  $ are polynomials of
degree at most $N\left(  g-1\right)  +N-1=Ng-1$.

We will use the notation $A\left(  z\right)  =\left(  A_{i,j}\left(  z\right)
\right)  _{i,j=0}^{N-1}$ for the loop-group element $A\colon\mathbb{T}%
\rightarrow\mathrm{U}\left(  N\right)  $. Since the Fourier expansion is
finite, $A\left(  z\right)  $ has the form%
\begin{equation}
A\left(  z\right)  =\sum_{k=0}^{g-1}z^{k}A^{\left(  k\right)  },
\label{eqRep.18}%
\end{equation}
where $A^{\left(  k\right)  }\in\mathcal{B}\left(  \mathbb{C}^{N}\right)  $
for $k=0,\dots,g-1$. The factorization in Lemma \ref{lemma2.3} motivates the
name \emph{genus} for $g$.

We have

\begin{corollary}
\label{CorRep.2}If $A\left(  z\right)  $ is a general polynomial of $z$ with
values in $\mathcal{B}\left(  \mathbb{C}^{N}\right)  $ of the form
\textup{(\ref{eqRep.18}),} the following four conditions
\textup{(\ref{eqRep.19})--(\ref{eqRep.22})} are equivalent:%
\begin{align}
&  \begin{minipage}[t]{\displayboxwidth}\raggedright$A\left( z\right) ^{\ast
}A\left( z\right) =\openone_{N}$, $z\in\mathbb{T}%
$, i.e., $A$ takes values in $\mathrm{U}\left( N\right) $; \end{minipage}%
\label{eqRep.19}\\
&  \begin{minipage}[t]{\displayboxwidth}\raggedright$\sum_{k}A^{\left
( k\right) \,\ast}A^{\left( k+n\right) }= \begin{cases} \openone_{N}
& \text{ if }n=0, \\ 0 & \text{ if }n\in\mathbb{Z}\setminus\left
\{ 0\right\} , \end{cases} $ with the convention that $A^{\left( m\right
) }=0$ if $m\notin\left\{ 0,1,\dots,g-1\right\} $; \end{minipage}%
\label{eqRep.20}\\
&  \begin{minipage}[t]{\displayboxwidth}\raggedright
there are projections $P_{1},\dots,P_{s}$ in $\mathcal{B}\left( \mathbb{C}%
^{N}\right) $, positive integers $r_{1},\dots,r_{s}$, and a unitary $W\in
\mathrm{U}\left( N\right) $ such that $A\left( z\right) =\left( \prod
_{j=1}^{s}\left( \openone_{N}-P_{j}+z^{r_{j}}P_{j}\right) \right
) W$; \end{minipage} \label{eqRep.21}%
\end{align}
and%
\begin{equation}
\begin{minipage}[t]{\displayboxwidth}\raggedright there are projections $Q_{1}%
,Q_{2},\dots,Q_{g-1}$ and a unitary $V\in\mathrm{U}\left( N\right
) $ such that \end{minipage} \label{eqRep.22}%
\end{equation}%
\begin{align*}
A^{\left(  0\right)  }  &  =V\prod_{j=1}^{g-1}\left(  \openone_{N}%
-Q_{j}\right)  ,\\
A^{\left(  1\right)  }  &  =V\sum_{j=1}^{g-1}\left(  \openone_{N}%
-Q_{1}\right)  \cdots\\
&  \qquad\cdots\left(  \openone_{N}-Q_{j-1}\right)  Q_{j}\left(  \openone
_{N}-Q_{j+1}\right)  \cdots\\
&  \qquad\qquad\cdots\left(  \openone_{N}-Q_{g-1}\right)  ,\\
\vdots &  \mathrel{\phantom{=W}}\vdots\\
\qquad A^{\left(  g-1\right)  }  &  =V\prod_{j=1}^{g-1}Q_{j}.
\end{align*}
\end{corollary}

\begin{proof}
This follows from Proposition \ref{proposition2.2} and explicit computations.
\end{proof}

\begin{remark}
\label{RemRep.3}The case $g=2=N$ includes the family of wavelets introduced by
Daubechies \cite{Dau92} and studied further in \cite{BEJ00}. Note that $g=2$
yields the representation%
\begin{equation}
A^{\left(  0\right)  }=V\left(  \openone_{N}-Q\right)  ,\qquad A^{\left(
1\right)  }=VQ, \label{eqRep.23}%
\end{equation}
by \textup{(\ref{eqRep.22}).} But then \textup{(\ref{eqRep.20})} takes the
form%
\begin{equation}
A^{\left(  0\right)  \,\ast}A^{\left(  0\right)  }=\openone_{N}-Q,\qquad
A^{\left(  1\right)  \,\ast}A^{\left(  1\right)  }=Q, \label{eqRep.24}%
\end{equation}
which will be used in the sequel.
\end{remark}

Let us return to the connection between loop-group elements and
representations of $\mathcal{O}_{N}$. The algebra $\mathcal{O}_{N}$ has the
following basic representation $s_{i}\mapsto S_{i}$ on $L^{2}\left(
\mathbb{T}\right)  $:%
\begin{equation}
S_{j}\xi\left(  z\right)  =z^{j}\xi\left(  z^{N}\right)  ,\qquad0\leq
j<N,\;\xi\in L^{2}\left(  \mathbb{T}\right)  ,\;z\in\mathbb{T}\text{.}
\label{eqRep.25}%
\end{equation}
If $s_{i}\mapsto T_{i}$ is any other representation, then as noted in
\cite{Cun77}, \cite{Cun80}, and \cite{BJP96}, the system $\left(  S_{j}^{\ast
}T_{i}^{{}}\right)  _{i,j}$ realizes a unitary $N\times N$ matrix with entries
in $\mathcal{B}\left(  L^{2}\left(  \mathbb{T}\right)  \right)  $.

\begin{proposition}
\label{ProRep.4}Let $\left(  S_{i}\right)  $ be the basic representation of
$\mathcal{O}_{N}$, and let $\left(  T_{i}\right)  $ be an arbitrary
representation. Then the following two conditions are equivalent:%
\begin{equation}
\begin{minipage}[t]{\displayboxwidth}\raggedright Each operator $S_{j}^{\ast
}T_{i}^{{}}$ on $L^{2}\left( \mathbb{T}\right
) $ is a multiplication operator; \end{minipage} \label{eqRep.26}%
\end{equation}
and%
\begin{equation}
\begin{minipage}[t]{\displayboxwidth}\raggedright the $T_{i}%
$-representation has the form $$ T_{i}\xi\left( z\right) =m_{i}\left
( z\right) \xi\left( z^{N}\right) ,\qquad0\leq i<N,\;\xi\in L^{2}%
\left( \mathbb{T}\right) ,\;z\in\mathbb{T}, $$ where $m_{0},\dots,m_{N-1}\in
L^{\infty}\left( \mathbb{T}\right) $. \end{minipage} \label{eqRep.27}%
\end{equation}
\end{proposition}

\begin{proof}
\underline{(\ref{eqRep.26}) $\Rightarrow$ (\ref{eqRep.27})}: If
(\ref{eqRep.26}) holds, then there are functions $A_{i,j}\in L^{\infty}\left(
\mathbb{T}\right)  $ such that $S_{j}^{\ast}T_{i}^{{}}=A_{i,j}$ where the
right-hand side also denotes the multiplication operator determined by the
function in question. Then%
\begin{equation}
T_{i}=\sum_{j=0}^{N-1}S_{j}^{{}}S_{j}^{\ast}T_{i}^{{}}=\sum_{j=0}^{N-1}%
S_{j}A_{i,j}=\sum_{j=0}^{N-1}A_{i,j}\left(  z^{N}\right)  S_{j}.
\label{eqRep.28}%
\end{equation}
Setting%
\begin{equation}
m_{i}\left(  z\right)  =\sum_{j=0}^{N-1}A_{i,j}\left(  z^{N}\right)  z^{j}
\label{eqRep.29}%
\end{equation}
and using the Cuntz relations, we conclude that $T_{i}$ has the form in
(\ref{eqRep.27}).

\underline{(\ref{eqRep.27}) $\Rightarrow$ (\ref{eqRep.26})}: If, conversely,
the representation $\left(  T_{i}\right)  $ is given to have the form in
(\ref{eqRep.27}), then one checks that the filter functions $m_{i}$ must
satisfy the unitarity property (\ref{eq1.22}) above. If then $A_{i,j}\left(
z\right)  $ are given by (\ref{eq1.23}), then%
\begin{equation}
S_{j}^{\ast}T_{i}^{{}}\xi\left(  z\right)  =\frac{1}{N}\sum_{w^{N}=z}%
w^{-j}T_{i}\xi\left(  w\right)  =\frac{1}{N}\sum_{w^{N}=z}w^{-j}m_{i}\left(
w\right)  \xi\left(  z\right)  , \label{eqRep.30}%
\end{equation}
and we conclude that (\ref{eqRep.26}) holds, i.e., $S_{j}^{\ast}T_{i}^{{}}$ is
multiplication by the function $A_{i,j}\left(  \,\cdot\,\right)  $, where
$A_{i,j}$ is derived from $m_{i}$ via (\ref{eq1.23})
\end{proof}

We will henceforth denote the representation $T_{i}$ defined by
(\ref{eqRep.30}) by $T_{i}^{\left(  A\right)  }$, so in particular,%
\begin{equation}
T_{i}^{\left(  \openone\right)  }=S_{i}^{{}}, \label{eqRep.31}%
\end{equation}
where the notation $S_{i}$ is reserved for the representation defined by
(\ref{eqRep.25}).

We will now apply the results in (\ref{eqRep.5})--(\ref{eqRep.17}) to analyze
these representations further. As mentioned after (\ref{eqRep.17}), the linear
span $\mathcal{K}$ of $1,z^{-1},z^{-2},\dots,z^{-r}$ is cyclic and
$T^{\ast}$-invariant
for a suitable $r\in\mathbb{N}$,
i.e., satisfies (\ref{eqRep.5})--(\ref{eqRep.6}) with
$T$ in lieu of $S$, see \cite[Proposition 3.1 and Corollary 3.3]{BEJ00}.
If $m_{0},\dots,m_{N-1}$ can
be derived from a polynomial loop $A\left(  z\right)  $ of degree $g-1$ by
(\ref{eq1.24}), we may
explicitly estimate $r$.
To this end, let us look at the action of
$T_{i}^{\left(  A\right)  \,\ast}$ on $e_{n}\left(  z\right)  =z^{n}$. Put
$n=j+Nl$ uniquely, where $j\in\left\{  0,1,\dots,N-1\right\}  $ and
$l\in\mathbb{Z}$. Then, using (\ref{eqRep.28}),%
\begin{align}
T_{i}^{\left(  A\right)  \,\ast}e_{n}  &  =\sum_{j^{\prime}}S_{j^{\prime}%
}^{\ast}\overline{A_{i,j^{\prime}}\left(  z^{N}\right)  }e_{n}
\label{eqRep.32}\\
&  =\sum_{k}\sum_{j^{\prime}}S_{j^{\prime}}^{\ast}\overline{A_{i,j^{\prime}%
}^{\left(  k\right)  }}z^{-Nk}e_{n}=\sum_{k}\sum_{j^{\prime}}\overline
{A_{i,j^{\prime}}^{\left(  k\right)  }}S_{j^{\prime}}^{\ast}e_{n-Nk}.\nonumber
\end{align}
But $n-Nk=j+Nl-Nk=j+N\left(  l-k\right)  $, and since%
\begin{equation}
S_{j^{\prime}}^{\ast}e_{j+N\left(  l-k\right)  }=%
\begin{cases}
e_{l-k} & \text{\quad if }j^{\prime}=j,\\
0 & \text{\quad if }j^{\prime}\neq j,
\end{cases}
\label{eqRep.33}%
\end{equation}
we obtain%
\begin{equation}
T_{i}^{\left(  A\right)  \,\ast}e_{j+Nl}=\sum_{k=0}^{g-1}\overline
{A_{i,j}^{\left(  k\right)  }}e_{l-k}. \label{eqRep.34}%
\end{equation}
It follows that $T_{i}^{\left(  A\right)  \,\ast}$ is represented by a matrix
which is a slanted block matrix of the following form:%
\begin{equation}%
\begin{array}
[c]{cc}
&
%TCIMACRO{\TeXButton{ShadedRectangle}{\savebox{\shaderect}(3,5)
%{\begin{picture}(3,5)\multiput(0,0.25)(0,1){3}{\line(6,5){3}}
%\multiput(0.9,0)(-0.9,3.25){2}{\line(6,5){2.1}}\multiput(2.1,0)(-2.1,4.25){2}
%{\line(6,5){0.9}}\end{picture}}}}%
%BeginExpansion
\savebox{\shaderect}(3,5)
{\begin{picture}(3,5)\multiput(0,0.25)(0,1){3}{\line(6,5){3}}
\multiput(0.9,0)(-0.9,3.25){2}{\line(6,5){2.1}}\multiput(2.1,0)(-2.1,4.25){2}
{\line(6,5){0.9}}\end{picture}}%
%EndExpansion
\begin{picture}(24,16)(-12,-10) \multiput(-12,-8)(4,1){6}{\frame
{\usebox{\shaderect}}} \put(-13,0){\vector(1,0){25}} \put(0,-10){\vector
(0,1){16}} \multiput(0,-0.033)(0.1,0.1){61}{\makebox(0,0){$\cdot$}}
\multiput(0,-0.033)(-0.1,-0.1){101}{\makebox(0,0){$\cdot$}} \put
(0,-0.028){\makebox(0,0){$\circ$}} \put(-7,-7.028){\makebox(0,0){$\circ$}}
\put(-0.1,0.15){\makebox(0,0)[br]{$\scriptstyle(0,0)$}} \put
(-6.9,-7.15){\makebox(0,0)[tl]{$\scriptstyle(\mkern-2mu -\mkern-2mu r,\mkern
-2mu -\mkern-2mu r)$}} \put(0.1,-5.15){\makebox(0,0)[tl]{$\scriptstyle
\mkern-2mu -\mkern-2mu (g\mkern-2mu -\mkern-2mu 1)$}} \put(3.06,0.15){\makebox
(0,0)[br]{$\scriptscriptstyle N\mkern-2mu -\mkern-2mu 1$}} \put
(4.13,0.1){\makebox(0,0)[bl]{$\scriptscriptstyle N$}} \put(8.16,0.1){\makebox
(0,0)[bl]{$\scriptscriptstyle2\mkern-2mu N$}} \end{picture}%
\end{array}
\label{eqRep.35}%
\end{equation}
The matrix elements outside the shaded blocks are all zero, and each block has
$N$ columns and $g$ rows. Each block is a translate $N$ steps to the right and
$1$ step up compared to the previous one,and the most central block is located
with corners at $\left(  0,0\right)  $, $\left(  N-1,0\right)  $, $\left(
N-1,-g+1\right)  $, $\left(  0,-g+1\right)  $. It follows that if we fix one
$n\in\mathbb{Z}$, iterated applications of $T_{i}^{\left(  A\right)  \,\ast}$
for various $i$'s ultimately will transform this vector into a linear
combination of $e_{m}$'s where $m\in\left\{  0,-1,-2,\dots,-r\right\}  $,
where $r$ is the largest integer such that $\left(  -r,-r\right)  $ is
contained in one of the blocks. Thus the linear span $\mathcal{K}$ of
$\left\{  e_{0},e_{-1},\dots,e_{-r}\right\}  $ is $T_{i}^{\left(  A\right)
\,\ast}$-invariant, and because of the relation
\begin{equation}
\openone=\sum_{\substack{I\\\left|  I\right|  =n}}T_{I}^{\left(  A\right)
}T_{I}^{\left(  A\right)  \,\ast}, \label{eqRep.35bis}%
\end{equation}
$\mathcal{K}$ will also be cyclic (see \cite[Section 3]{BEJ00} for more
details of this argument). Thus, computing $r$ more explicitly, we deduce

\begin{proposition}
\label{ProRep.5}If $T^{\left(  A\right)  }$ is the representation of the Cuntz
algebra $\mathcal{O}_{N}$ defined by a polynomial loop $A\left(  z\right)  $
of degree $g-1$, then the subspace%
\begin{equation}
\mathcal{K}=\operatorname*{lin}\operatorname*{span}\left\{  e_{0}%
,e_{-1},e_{-2},\dots,e_{-r}\right\}  \label{eqRep.36}%
\end{equation}
satisfies \textup{(\ref{eqRep.5})} and \textup{(\ref{eqRep.6}),} i.e.,
$\mathcal{K}$ is $T^{\left(  A\right)  \,\ast}$-invariant and cyclic, where%
\begin{equation}
r=g+\left\lfloor \frac{g-1}{N-1}\right\rfloor 
=\left\lfloor \frac{gN-1}{N-1}\right\rfloor . \label{eqRep.37}%
\end{equation}
Here $\left\lfloor x\right\rfloor $ is the largest integer $\leq x$.
\end{proposition}

\begin{proof}
{}From the figure (\ref{eqRep.35}) it follows that%
\begin{equation}
r=%
\left\{
\begin{array}
[c]{lll}%
g & \text{if }g<N, & \text{i.e., }g-1<N-1,\\
g+1 & \text{if }N+1\leq g+1<2N, & \text{i.e., }N-1\leq g-1<2N-2,\\
g+2 & \text{if }2N+1\leq g+2<3N, & \text{i.e., }2N-2\leq g-1<3N-3,\\
\text{etc.} &  &
\end{array}
\right. 
\label{eqRep.37bis}%
\end{equation}
Here the first column of ranges (``if $\dots$'')
is taken from the box diagram (\ref{eqRep.35}),
and the second column of ranges (``i.e., $\dots$'')
is
derived from the first by subtraction of suitable integers.
The second column
shows how the ratio $\frac{g-1}{N-1}$ arises in (\ref{eqRep.37}).
\end{proof}

\begin{remark}
\label{RemRepNew.6}Proposition \textup{\ref{ProRep.5}} could
alternatively have been proved by
exactly the method used to show
\textup{(4.18)} in \cite{BEJ00} from Remark \textup{3.2}
there, using joint invariant sets
for the maps $n\rightarrow \left( n-k\right) /N$ for
$k=0,1,\dots ,gN-1$. Note that a
consequence of \textup{(\ref{eq1.23})} and \textup{(\ref{eq1.24})}
is that the degree of all the polynomials
$A_{i,j}\left( z\right) $ for $i,j=0,\dots ,N-1$ is
at most $g-1$ if and only if
the degree of all the polynomials
$m_{i}\left( z\right) $ for $i=0,\dots ,N-1$ is at
most $gN-1$.
\end{remark}

The case considered in detail in \cite{BEJ00} was $N=g=2$, so we regain the
result $r=3$ from there, i.e., in that case $\mathcal{K}=\operatorname*{span}%
\left\{  e_{0},e_{-1},e_{-2},e_{-3}\right\}  $.
Note that the denominator $N-1$ in (\ref{eqRep.37})
shows the dependence of $r=\left( \dim\mathcal{K}\right) -1$
on the scaling number $N$
for the given wavelet filter. So
when $N>2$, the denominator
$N-1$ in (\ref{eqRep.37}) yields a smaller value for
$r$ than the value $gN-1$
from \cite{BEJ00} in the special case $N=2$.

We will now apply Theorem \ref{ThmRep.1} to the representations above, when
$\mathcal{K}=\mathcal{K}^{\left(  A\right)  }$ is given by (\ref{eqRep.36})
and $r=r^{\left(  A\right)  }$ by (\ref{eqRep.37}). The operators $V_{i}%
\in\mathcal{B}\left(  \mathcal{K}^{\left(  A\right)  }\right)  $ defined by
$T_{i}^{\left(  A\right)  }$ in lieu of $S_{i}$ in (\ref{eqRep.7}) will be
denoted by $V_{i}^{\left(  A\right)  }$, and then the $\sigma$ defined by
(\ref{eqRep.9}) is denoted by $\sigma^{\left(  A\right)  }=\sigma^{\left(
A,A\right)  }$, i.e.,
\begin{equation}
\sigma^{\left(  A,A\right)  }\left(  X\right)  =\sum_{i=0}^{N-1}V_{i}^{\left(
A\right)  }XV_{i}^{\left(  A\right)  \,\ast}\label{eqRep.38}%
\end{equation}
for $X\in\mathcal{B}\left(  \mathcal{K}^{\left(  A\right)  }\right)  $. If
$T^{\left(  A\right)  }$, $T^{\left(  B\right)  }$ are two representations of
our kind, let $\sigma^{\left(  B,A\right)  }$ denote the $\rho$ defined by
(\ref{eqRep.14}), i.e.,%
\begin{equation}
\sigma^{\left(  B,A\right)  }\left(  X\right)  =\sum_{i}V_{i}^{\left(
B\right)  }XV_{i}^{\left(  A\right)  \,\ast},\label{eqRep.39}%
\end{equation}
so $\sigma^{\left(  B,A\right)  }$ is in $\mathcal{B}\left(  \mathcal{K}%
^{\left(  A\right)  },\mathcal{K}^{\left(  B\right)  }\right)  $. We will
mostly be interested in the situation that $A$ and $B$ have the same genus
$g$, and then $\mathcal{K}^{\left(  A\right)  }=\mathcal{K}^{\left(  B\right)
}$ by Proposition \ref{ProRep.5} (of course we may always replace $g^{\left(
A\right)  }$, $g^{\left(  B\right)  }$ by $\max\left\{  g^{\left(  A\right)
},g^{\left(  B\right)  }\right\}  $ and thus assume $\mathcal{K}^{\left(
A\right)  }=\mathcal{K}^{\left(  B\right)  }$). Now $\mathcal{B}\left(
\mathcal{K}^{\left(  A\right)  },\mathcal{K}^{\left(  B\right)  }\right)  $
can be made into a Hilbert space $\mathcal{K}^{\left(  B,A\right)  }$ in a
natural fashion by defining the norm of $\psi\in\mathcal{B}\left(
\mathcal{K}^{\left(  A\right)  },\mathcal{K}^{\left(  B\right)  }\right)  $ by%
\begin{equation}
\left\|  \psi\right\|  _{2}^{\left(  B,A\right)  }=\operatorname*{Tr}%
\nolimits^{\left(  A\right)  }\left(  \psi^{\ast}\psi\right)
=\operatorname*{Tr}\nolimits^{\left(  B\right)  }\left(  \psi\psi^{\ast
}\right)  ,\label{eqRep.40}%
\end{equation}
where $\operatorname*{Tr}^{\left(  A\right)  }$, $\operatorname*{Tr}^{\left(
B\right)  }$ are the standard unnormalized traces on $\mathcal{B}\left(
\mathcal{K}^{\left(  A\right)  }\right)  $, $\mathcal{B}\left(  \mathcal{K}%
^{\left(  B\right)  }\right)  $, respectively. Then $\sigma^{\left(
B,A\right)  }\in\mathcal{B}\left(  \mathcal{K}^{\left(  B,A\right)  }\right)
$. Because of Theorem \ref{ThmRep.1}, our main concern will be the spectral
multiplicity of $1$ in the spectrum $\operatorname*{Sp}\left(  \sigma^{\left(
B,A\right)  }\right)  $, as well as the associated eigensubspace
$\mathcal{K}^{\left(  B,A\right)  }$. It turns out that these items are much
simpler to compute for the adjoint $\sigma^{\left(  B,A\right)  \,\ast}%
\in\mathcal{B}\left(  \mathcal{K}^{\left(  B,A\right)  }\right)  $ (i.e., the
adjoint of $\sigma^{\left(  B,A\right)  }$ as a Hilbert-space operator). Then
$\operatorname*{Sp}\left(  \sigma^{\left(  B,A\right)  \,\ast}\right)
=\overline{\operatorname*{Sp}\left(  \sigma^{\left(  B,A\right)  }\right)  }$,
and in particular $1\in\operatorname*{Sp}\left(  \sigma^{\left(  B,A\right)
}\right)  $ if and only if $1\in\operatorname*{Sp}\left(  \sigma^{\left(
B,A\right)  \,\ast}\right)  $. If $A=B$ and the completely positive unital map
$\sigma^{\left(  A,A\right)  }$ admits a faithful invariant state, it follows
from \cite[Lemma 6.3]{BrJo97a} that the dimensions of the fixed-point sets of
$\sigma^{\left(  A,A\right)  }$ and $\sigma^{\left(  A,A\right)  \,\ast}$ are
the same. However, we will see that $\sigma^{\left(  A,A\right)  }$ does not
necessarily admit a faithful invariant state. Our strategy will be to compute
$\operatorname*{Sp}\left(  \sigma^{\left(  B,A\right)  }\right)  $ by
evaluating%
\begin{equation}
\det\left(  x\openone-\sigma^{\left(  A,B\right)  \,\ast}\right)
=\overline{\det\left(  \bar{x}\openone-\sigma^{\left(  A,B\right)  }\right)
}\label{eqRep.41}%
\end{equation}
and then compute the fixed points by hand. Note that from (\ref{eqRep.39}) we
have%
\begin{equation}
\sigma^{\left(  B,A\right)  \,\ast}\left(  X\right)  =\sum_{i}V_{i}^{\left(
B\right)  \,\ast}XV_{i}^{\left(  A\right)  }\label{eqRep.42}%
\end{equation}
for $X\in\mathcal{K}^{\left(  B,A\right)  }=\mathcal{B}\left(  \mathcal{K}%
^{\left(  A\right)  },\mathcal{K}^{\left(  B\right)  }\right)  $. If
$X=E_{-k,-l}=$ the rank-one partial isometry mapping $e_{-l}$ into $e_{-k}$,
we compute, in Dirac's notation,%
\begin{equation}
\sigma^{\left(  B,A\right)  }\left(  E_{-k,-l}\right)  =\sum_{i}V_{i}^{\left(
B\right)  }E_{-k,-l}^{{}}V_{i}^{\left(  A\right)  \,\ast}=\sum_{i}\left|
V_{i}^{\left(  B\right)  }e_{-k}^{{}}\right\rangle \left\langle V_{i}^{\left(
A\right)  }e_{-l}^{{}}\right|  ,\label{eqRep.43}%
\end{equation}
which is complicated to handle. On the other hand,
\begin{multline}
\sigma^{\left(  B,A\right)  \,\ast}\left(  E_{-k,-l}\right)   =\sum
_{i}V_{i}^{\left(  B\right)  \,\ast}E_{-k,-l}V_{i}^{\left(  A\right)
}\label{eqRep.44}\\
=\sum_{i}\left|  V_{i}^{\left(  B\right)  \,\ast}e_{-k}\right\rangle
\left\langle V_{i}^{\left(  A\right)  \,\ast}e_{-l}\right|  
=\sum_{i}\left|  T_{i}^{\left(  B\right)  \,\ast}e_{-k}\right\rangle
\left\langle T_{i}^{\left(  A\right)  \,\ast}e_{-l}\right|  ,
\end{multline}
where the last step used $T^{\ast}$-invariance of $\mathcal{K}$, see
(\ref{eqRep.5}) and (\ref{eqRep.7}). But $T^{\ast}$ is given by
(\ref{eqRep.34}), and hence the matrix elements of $\sigma^{\left(
B,A\right)  }$ can be computed explicitly. Let us do the calculation in a
special case:

\subsection{The case $g=2$ and $N>g$}

Then $r=2$ by Proposition \ref{ProRep.5}, and $\mathcal{K}$ is
three-dimensional. Then the loop $A\left(  z\right)  $ has the form%
\begin{equation}
A\left(  z\right)  =V\left(  1-Q+zQ\right)  =A^{\left(  0\right)
}+zA^{\left(  1\right)  }, \label{eqRep.45}%
\end{equation}
where $V$ is a unitary and $Q$ a projection in $\mathcal{B}\left(
\mathbb{C}^{N}\right)  $. As noted in (\ref{eqRep.24}), we have%
\begin{equation}
A^{\left(  0\right)  \,\ast}A^{\left(  0\right)  }=\openone_{N}-Q,\qquad
A^{\left(  1\right)  \,\ast}A^{\left(  1\right)  }=Q. \label{eqRep.46}%
\end{equation}
Now define%
\begin{equation}
\lambda_{i,j}=\left(  \openone-Q\right)  _{i,j}=\delta_{i,j}-Q_{i,j},\qquad
i,j=0,1,2,\dots,N-1. \label{eqRep.47}%
\end{equation}
Setting $\lambda_{i}:=\lambda_{i,i}$, then (\ref{eqRep.44}) and
(\ref{eqRep.34}) give%
\begin{equation}
\left\{
\begin{array}
[c]{l}%
\sigma^{\left(  A,A\right)  \,\ast}\left(  E_{0,0}\right)  =\lambda_{0}%
E_{0,0}+\left(  1-\lambda_{0}\right)  E_{-1,-1},\\
\vphantom{\vdots}\sigma^{\left(  A,A\right)  \,\ast}\left(  E_{-1,-1}\right)
=\lambda_{N-1}E_{-1,-1}+\left(  1-\lambda_{N-1}\right)  E_{-2,-2},\\
\vphantom{\vdots}\sigma^{\left(  A,A\right)  \,\ast}\left(  E_{-2,-2}\right)
=\lambda_{N-2}E_{-1,-1}+\left(  1-\lambda_{N-2}\right)  E_{-2,-2},\\
\vphantom{\vdots}\sigma^{\left(  A,A\right)  \,\ast}\left(  E_{0,-1}\right)
=\lambda_{0,N-1}E_{0,-1}-\lambda_{0,N-1}E_{-1,-2},\\
\vphantom{\vdots}\sigma^{\left(  A,A\right)  \,\ast}\left(  E_{0,-2}\right)
=\lambda_{0,N-2}E_{0,-1}-\lambda_{0,N-2}E_{-1,-2},\\
\vphantom{\vdots}\sigma^{\left(  A,A\right)  \,\ast}\left(  E_{-1,-2}\right)
=\lambda_{N-1,N-2}E_{-1,-1}-\lambda_{N-1,N-2}E_{-2,-2},\\
\vphantom{\vdots}\sigma^{\left(  A,A\right)  \,\ast}\left(  E_{-k,-l}\right)
=\sigma^{\left(  A,A\right)  \,\ast}\left(  E_{-l,-k}\right)  ^{\ast
}\text{\qquad whenever }0\leq l<k\leq2.
\end{array}
\right.  \label{eqRep.48}%
\end{equation}
{}From here one computes%
\begin{multline}
\det\left(  x\openone-\sigma^{\left(  A,A\right)  }\right)
=\det\left(  x\openone-\sigma^{\left(  A,A\right)  \,\ast}\right)
\label{eqRep.49}\\
=x^{4}\left(  x-1\right)  \left(  x-\lambda_{0}\right)  \left(
x-\left(  \lambda_{N-1}-\lambda_{N-2}\right)  \right)  \left(  x-\lambda
_{0,N-1}\right)  \left(  x-\overline{\lambda_{0,N-1}}\right)  .
\end{multline}
Since the $\lambda_{i,j}$ are matrix elements of a projection, we see
immediately that the eigenvalue $1$ has multiplicity one except in special
cases as when $\lambda_{0}=1$, $\lambda_{N-1}=1$ and $\lambda_{N-2}=0$, or
$\lambda_{0,N-1}=1$. Hence:

\begin{theorem}
\label{ThmRep.6}If $g=2$ and $N>2$, the representation $T^{\left(  A\right)
}$ is irreducible for generic loops $A$. It is reducible if and only if the
projection $Q\in\mathcal{B}\left(  \mathbb{C}^{N}\right)  $ defining $A$ by
\textup{(\ref{eqRep.45})} satisfies at least one of the conditions%
\begin{equation}
Q_{0,0}=0\qquad(\iff\lambda_{0}=1) \label{eqRep.50}%
\end{equation}
or%
\begin{equation}
Q_{N-1,N-1}=0\text{\quad and\quad}Q_{N-2,N-2}=1\qquad(\iff\lambda
_{N-1}=1\text{ and }\lambda_{N-2}=0). \mkern-12mu\label{eqRep.51}%
\end{equation}
\end{theorem}

\begin{proof}
We will use the fact that the $\lambda_{i,j}$ numbers are the matrix entries
of a projection, $\openone-Q$. Since $\left(  \openone-Q\right)  ^{2}=\left(
\openone-Q\right)  $, we have%
\begin{equation}
\lambda_{k}=\lambda_{k}^{2}+\sum_{j\neq k}\left|  \lambda_{k,j}\right|  ^{2}.
\label{eqRep.53}%
\end{equation}
Thus the case (\ref{eqRep.50}) may occur with%
\begin{equation}
Q=\left(
\begin{tabular}
[c]{c|c}%
$0$ & $%
\begin{array}
[c]{cccc}%
0 & 0 & \cdots & 0
\end{array}
$\\\hline
$%
\begin{array}
[c]{c}%
0\\
0\\
\vdots\\
0
\end{array}
$ & $\rule{15pt}{0pt}%
\begin{array}
[c]{ccc}%
\hbox to0pt{\hss$\,\cdots$\hss}\hbox to0pt{\hss\raisebox{-2pt}{$\vdots$}\hss}
& \cdots & \hbox to0pt{\hss$\,\cdots$\hss}\hbox to0pt{\hss\raisebox
{-2pt}{$\vdots$}\hss}\\
\vdots & Q_{\operatorname*{red}} & \vdots\\
\hbox to0pt{\hss$\,\cdots$\hss}\hbox to0pt{\hss\raisebox{-2pt}{$\vdots$}\hss}
& \cdots & \hbox to0pt{\hss$\,\cdots$\hss}\hbox to0pt{\hss\raisebox
{-2pt}{$\vdots$}\hss}%
\end{array}
\rule[-33pt]{0pt}{66pt}\rule{15pt}{0pt}$%
\end{tabular}
\right)  =\left(  0\right)  \oplus Q_{\operatorname*{red}}.
\label{eqRep.53bis}%
\end{equation}
Also (\ref{eqRep.51}) may occur with%
\begin{equation}
Q=\left(
\begin{tabular}
[c]{c|c}%
$%
\begin{array}
[c]{ccc}%
\hbox to0pt{\hss$\,\cdots$\hss}\hbox to0pt{\hss\raisebox{-2pt}{$\vdots$}\hss}
& \cdots & \hbox to0pt{\hss$\,\cdots$\hss}\hbox to0pt{\hss\raisebox
{-2pt}{$\vdots$}\hss}\\
\vdots & Q_{\operatorname*{red}} & \vdots\\
\hbox to0pt{\hss$\,\cdots$\hss}\hbox to0pt{\hss\raisebox{-2pt}{$\vdots$}\hss}
& \cdots & \hbox to0pt{\hss$\,\cdots$\hss}\hbox to0pt{\hss\raisebox
{-2pt}{$\vdots$}\hss}%
\end{array}
$ & $%
\begin{array}
[c]{cc}%
0 & 0\\
0 & 0\\
\vdots & \vdots\\
0 & 0
\end{array}
$\\\hline
$%
\begin{array}
[c]{cccc}%
0 & 0 & \cdots & 0\\
0 & 0 & \cdots & 0
\end{array}
$ & $%
\begin{array}
[c]{cc}%
1 & 0\\
0 & 0
\end{array}
$%
\end{tabular}
\right)  =Q_{\operatorname*{red}}\oplus\left(
\begin{array}
[c]{cc}%
1 & 0\\
0 & 0
\end{array}
\right)  . \label{eqRep.53ter}%
\end{equation}
But the case $Q_{0,N-1}=-1$, i.e., $\lambda_{0,N-1}=\lambda_{N-1,0}=1$, cannot
occur.
To see this, recall that, by (\ref{eqRep.24}) and (\ref{eqRep.47}),
the $\lambda_{i,j}$ numbers are the matrix
entries of $\openone_{N}-Q$, which is a
projection in $\mathbb{C}^{N}$, since $Q$ is.
Specifically,
\begin{equation}
\lambda_{0,N-1}=-\ip{\varepsilon_{0}}{Q\varepsilon_{N-1}},
\label{eqRepFeb4.lambdaip}
\end{equation}
where $\left\{ \varepsilon_{i}\right\}_{i=0}^{N-1}$
is the canonical basis for $\mathbb{C}^{N}$.
So we are considering the possibility
of having an off-diagonal entry
$\lambda_{0,N-1}$ assume the value $1$, and
that is impossible. In fact, let
$P=P^{\ast}=P^{2}$ be any projection, and let
$\varepsilon $, $\varepsilon^{\prime}$
be given orthogonal unit vectors. Then
\begin{equation}
\left| \ip{\varepsilon}{P\varepsilon^{\prime}}\right| \leq\frac{1}{2},
\label{eqRepFeb4.pound}
\end{equation}
and so in particular the value $1$ is excluded.
Note that (\ref{eqRepFeb4.pound}) is clearly sharp: take
$
\begin{pmatrix}
1/2 & 1/2 \\
1/2 & 1/2
\end{pmatrix}
$.
To prove (\ref{eqRepFeb4.pound}), note that
\begin{equation}
\left| \ip{\varepsilon}{P\varepsilon^{\prime}}\right| ^{2}
+ \left| \ip{\varepsilon^{\prime}}{P\varepsilon^{\prime}}\right| ^{2}
\leq\left\| P\varepsilon^{\prime}\right\| ^{2},
\label{eqRepFeb4.normsquare}
\end{equation}
by Bessel. Since
$\ip{\varepsilon^{\prime}}{P\varepsilon^{\prime}}
=\left\| P\varepsilon^{\prime}\right\| ^{2}$, we get
\begin{multline}
\left| \ip{\varepsilon}{P\varepsilon^{\prime}}\right| 
\leq\sqrt{\left\| P\varepsilon^{\prime}\right\| ^{2}
-\left\| P\varepsilon^{\prime}\right\| ^{4}}
=\left\| P\varepsilon^{\prime}\right\| \cdot 
\left( 1-\left\| P\varepsilon^{\prime}\right\| ^{2}\right) ^{1/2} \label{eqRepFeb4.normbound} \\
=\left\| P\varepsilon^{\prime}\right\| \cdot 
\left\| \left( \openone -P\right) \varepsilon^{\prime}\right\|
\leq\frac{1}{2}\left( \left\| P\varepsilon^{\prime}\right\| ^{2}
+\left\| \left( \openone -P\right) \varepsilon^{\prime}\right\| ^{2} 
\right) =\frac{1}{2}\left\| \varepsilon^{\prime}\right\| ^{2}
=\frac{1}{2}.
\end{multline}
This proves (\ref{eqRepFeb4.pound}).
\end{proof}

\begin{remark}
\label{RemRepNew.8}So in summary, we may
define a transformation $\sigma$ for
any system of numbers $\lambda_{i,j}$
via \textup{(\ref{eqRep.48})}, or the matrix
\textup{(\ref{eqRep.56})} below, and its spectrum
will be given by \textup{(\ref{eqRep.49})}. We
are then interested in when points
from the spectrum of $\sigma$ can attain
the value $1$, and we saw
that some can, and others
cannot, viz.,
if the $\lambda_{i,j}$ are derived from
an $A$ with $g=2$, then we
get the following estimates on the
points $\lambda_{0}$, $\lambda_{N-1}-\lambda_{N-2}$, $\lambda_{0,N-1}$,
and $\lambda_{N-1,0}$ in the spectrum of
$\sigma$:
\[
0\leq \lambda_{0}\leq 1,
\qquad
-1\leq \lambda_{N-1}-\lambda_{N-2}\leq 1,
\] 
and,
by the middle step in \textup{(\ref{eqRepFeb4.normbound}),}
\[
\left| \lambda_{0,N-1}\right| 
=\left| \lambda_{N-1,0}\right| 
\leq\min\left\{ 
\left( \lambda_{0}\cdot\left( 1-\lambda_{0}\right) \right) ^{1/2},
\left( \lambda_{N-1}\cdot\left( 1-\lambda_{N-1}\right) \right) ^{1/2}\right\} 
\leq\frac{1}{2}.
\]
In particular, the spectral
radius of such a $\sigma =\sigma^{\left( A\right) }$
is at most $1$. Nonetheless,
we know that $\sigma^{\left( A\right) }$ is generally not
contractive in the Hilbert space
$\mathcal{B}\left( \mathcal{K}^{\left( A\right) }\right) $
of Hilbert-Schmidt
operators, although $\sigma^{\left( A\right) }$ obviously is contractive as a
completely positive map
on the $C^*$-algebra
$\mathcal{B}\left( \mathcal{K}^{\left( A\right) }\right) $.
\end{remark}

\begin{remark}
\label{RemRep.7}Hence for a given unitary $V\in\mathrm{U}\left(  N\right)  $
and projection $Q\in\mathcal{B}\left(  \mathbb{C}^{N}\right)  $, we see that
the eigenvalue $1$ may have the following multiplicities in the characteristic
polynomial for various values of $A\sim\left(  V,Q\right)  $:\medskip

\noindent$1$ \textbf{has multiplicity} $3$: This occurs in one
case.

\begin{case}
\label{CasRemRep.7mult3(1)}$\lambda_{0}=1$, $\lambda_{N-2}=0$, and
$\lambda_{N-1}=1$. The projection $\openone-Q$ then has the form%
\begin{equation}
\openone-Q=\left(
\begin{tabular}
[c]{c|c|c}%
$1$ & $%
\begin{array}
[c]{ccc}%
0 & \cdots & 0
\end{array}
$ & $%
\begin{array}
[c]{cc}%
0 & 0
\end{array}
$\\\hline
$%
\begin{array}
[c]{c}%
0\\
\vdots\\
0
\end{array}
$ & $%
\begin{array}
[c]{ccc}%
\hbox to0pt{\hss$\,\cdots$\hss}\hbox to0pt{\hss\raisebox{-2pt}{$\vdots$}\hss}
& \cdots & \hbox to0pt{\hss$\,\cdots$\hss}\hbox to0pt{\hss\raisebox
{-2pt}{$\vdots$}\hss}\\
\vdots & P & \vdots\\
\hbox to0pt{\hss$\,\cdots$\hss}\hbox to0pt{\hss\raisebox{-2pt}{$\vdots$}\hss}
& \cdots & \hbox to0pt{\hss$\,\cdots$\hss}\hbox to0pt{\hss\raisebox
{-2pt}{$\vdots$}\hss}%
\end{array}
$ & $%
\begin{array}
[c]{cc}%
0 & 0\\
\vdots & \vdots\\
0 & 0
\end{array}
$\\\hline
$%
\begin{array}
[c]{c}%
0\\
0
\end{array}
$ & $%
\begin{array}
[c]{ccc}%
0 & \cdots & 0\\
0 & \cdots & 0
\end{array}
$ & $%
\begin{array}
[c]{cc}%
0 & 0\\
0 & 1
\end{array}
$%
\end{tabular}
\right)  , \label{eqRep.54}%
\end{equation}
where $P$ is a projection in $\mathcal{B}\left(  \mathbb{C}^{N-3}\right)  $
\textup{(}so when $N=3$ there is no more choice.\/\textup{)}\medskip
\end{case}

\noindent$1$ \textbf{has multiplicity} $2$: This occurs in two mutually
exclusive cases.

\begin{case}
\label{CasRemRep.7mult2(1)}$\lambda_{0}=1$ and \textup{(}$\lambda_{N-2}>0$ or
$\lambda_{N-1}<1$\textup{).}
\end{case}

\begin{case}
\label{CasRemRep.7mult2(2)}$\lambda_{0}<1$, $\lambda_{N-2}=0$, and
$\lambda_{N-1}=1$.
\end{case}

\noindent$1$ \textbf{has multiplicity} $1$: This occurs in all remaining cases.
\end{remark}

Since the dimension of the fixed-point set of $\sigma^{\left(  A\right)  }$ is
always at most equal to the multiplicity of $1$ in the characteristic
polynomial, it follows from Remark \ref{RemRep.7} and Theorem \ref{ThmRep.1}
that the linear dimension of the commutant of the $T^{\left(  A\right)  }%
$($=T^{\left(  V,Q\right)  }$)-representation is always at most $3$. Since the
smallest-dimensional nonabelian $C^{\ast}$-algebra is $\mathcal{B}\left(
\mathbb{C}^{2}\right)  $, which has dimension $4$, it follows that the
commutant is always abelian when $g=2$, and the representation decomposes into
a sum of at most three irreducible mutually disjoint representations. In order
to get the exact number, we have to compute the fixed-point space of
$\sigma^{\left(  A\right)  }$ exactly. To this end, we note that the matrix of
$\sigma^{\left(  A\right)  }$, relative to the basis $\left\{  E_{i,k};-2\leq
i,k\leq0\right\}  $ of $\mathcal{B}\left(  \mathbb{C}^{3}\right)  $ is (we
actually compute the matrix of $\sigma^{\left(  A\right)  \,\ast}$ using
(\ref{eqRep.48}) and take the adjoint)%
\begin{equation}
\addtolength{\tabcolsep}{-0.5\tabcolsep}\settowidth{\pairwidth}{$\scriptstyle
\left( 0,-2\right)  $}%
\begin{tabular}
[c]{r|ccc|ccc|ccc}%
$\scriptstyle\sigma^{\left(  A\right)  }$ & $\scriptstyle\left(  0,0\right)  $
& $\scriptstyle\left(  -1,-1\right)  $ & $\scriptstyle\left(  -2,-2\right)  $
& $\scriptstyle\left(  0,-1\right)  $ & $\scriptstyle\left(  0,-2\right)  $ &
$\scriptstyle\left(  -1,-2\right)  $ & $\scriptstyle\left(  -2,-1\right)  $ &
$\scriptstyle\left(  -1,0\right)  $ & $\scriptstyle\left(  -2,0\right)
$\\\hline
$\scriptstyle\left(  0,0\right)  $ & \rule{0pt}{14pt}\framebox{$\scriptstyle
\lambda_{0}$} & $\scriptstyle1-\lambda_{0}$ & $\scriptstyle0$ & $\scriptstyle
0$ & $\scriptstyle0$ & $\scriptstyle0$ & $\scriptstyle0$ & $\scriptstyle0$ &
$\scriptstyle0$\\
$\scriptstyle\left(  -1,-1\right)  $ & $\scriptstyle0$ & $\scriptstyle
\lambda_{N-1}$ & $\scriptstyle1-\lambda_{N-1}$ & $\scriptstyle0$ &
$\scriptstyle0$ & $\scriptstyle0$ & $\scriptstyle0$ & $\scriptstyle0$ &
$\scriptstyle0$\\
$\scriptstyle\left(  -2,-2\right)  $ & $\scriptstyle0$ & $\scriptstyle
\lambda_{N-2}$ & $\scriptstyle1-\lambda_{N-2}$ & $\scriptstyle0$ &
$\scriptstyle0$ & $\scriptstyle0$ & $\scriptstyle0$ & $\scriptstyle0$ &
$\scriptstyle0$\\\hline
$\scriptstyle\left(  0,-1\right)  $ & $\scriptstyle0$ & $\scriptstyle0$ &
$\scriptstyle0$ & \rule{0pt}{14pt}\rlap{\smash{\framebox{\begin{tabular}
[t]{cc}
$\scriptstyle\lambda_{N-1,0}$ & \makebox[\pairwidth]{$\scriptstyle0$}\\
$\scriptstyle\lambda_{N-2,0}$ & \makebox[\pairwidth]{$\scriptstyle0$}
\end{tabular}\hskip-3pt}}}\phantom{\hskip6pt$\scriptstyle\lambda_{N-1,0}$} &
\phantom{$\scriptstyle0$\hskip6pt} & $\scriptstyle-\lambda_{N-1,0}$ &
$\scriptstyle0$ & $\scriptstyle0$ & $\scriptstyle0$\\
$\scriptstyle\left(  0,-2\right)  $ & $\scriptstyle0$ & $\scriptstyle0$ &
$\scriptstyle0$ & \rule[-7pt]{0pt}{14pt}\phantom{\hskip6pt$\scriptstyle
\lambda_{N-2,0}$} & \phantom{$\scriptstyle0$\hskip6pt} & $\scriptstyle
-\lambda_{N-2,0}$ & $\scriptstyle0$ & $\scriptstyle0$ & $\scriptstyle0$\\
$\scriptstyle\left(  -1,-2\right)  $ & $\scriptstyle0$ & $\scriptstyle
\lambda_{N-2,N-1}$ & $\scriptstyle-\lambda_{N-2,N-1}$ & $\scriptstyle0$ &
$\scriptstyle0$ & $\scriptstyle0$ & $\scriptstyle0$ & $\scriptstyle0$ &
$\scriptstyle0$\\\hline
$\scriptstyle\left(  -2,-1\right)  $ & $\scriptstyle0$ & $\scriptstyle
\lambda_{N-1,N-2}$ & $\scriptstyle-\lambda_{N-1,N-2}$ & $\scriptstyle0$ &
$\scriptstyle0$ & $\scriptstyle0$ & $\scriptstyle0$ & $\scriptstyle0$ &
$\scriptstyle0$\\\cline{9-10}%
$\scriptstyle\left(  -1,0\right)  $ & $\scriptstyle0$ & $\scriptstyle0$ &
$\scriptstyle0$ & $\scriptstyle0$ & $\scriptstyle0$ & $\scriptstyle0$ &
$\scriptstyle-\lambda_{0,N-1}$ & \multicolumn{1}{|c}{$\scriptstyle
\lambda_{0,N-1}$} & \multicolumn{1}{c|}{$\scriptstyle0$}\\
$\scriptstyle\left(  -2,0\right)  $ & $\scriptstyle0$ & $\scriptstyle0$ &
$\scriptstyle0$ & $\scriptstyle0$ & $\scriptstyle0$ & $\scriptstyle0$ &
$\scriptstyle-\lambda_{0,N-2}$ & \multicolumn{1}{|c}{$\scriptstyle
\lambda_{0,N-2}$} & \multicolumn{1}{c|}{$\scriptstyle0$}\\\cline{9-10}%
\end{tabular}
\label{eqRep.56}%
\end{equation}
Thus one may compute the eigenvectors of $\sigma^{\left(  A\right)  }$. The
generic result is as in Table \ref{TabEigMult}.

\begin{table}[ptbh]
\caption{Eigenvalues, multiplicities, and eigenvectors of $\sigma^{\left(
A\right)  }$}%
\label{TabEigMult}
\begin{tabular}
[c]{l|l|l}%
Multiplicity & Eigenvalue & Eigenvectors\\\hline
$\vphantom{\vdots}4$ & $0$ & $E_{0,-2}$, $E_{-2,0}$, $E_{0,-1}+E_{-1,-2}$,
$E_{-1,0}+E_{-2,-1}$\\\hline
$\vphantom{\vdots}1$ & $\lambda_{0}$ & $E_{0,0}$\\\hline
$\vphantom{\vdots}1$ & $1$ & $\openone_{\mathcal{K}}=E_{0,0}+E_{-1,-1}%
+E_{-2,-2}$\\\hline
$\vphantom{\vdots}1$ & $\lambda_{N-1}-\lambda_{N-2}$ & \hskip-\tabcolsep$%
\begin{array}
[t]{l}%
\left(  1-\lambda_{0}\right)  \left(  1-\lambda_{N-1}\right)  E_{0,0}\\
\quad+\left(  \lambda_{N-1}-\lambda_{N-2}-\lambda_{0}\right)  \left(
1-\lambda_{N-1}\right)  E_{-1,-1}\\
\qquad-\lambda_{N-2}\left(  \lambda_{N-1}-\lambda_{N-2}-\lambda_{0}\right)
E_{-2,-2}%
\end{array}
$\\\hline
$\vphantom{\vdots}1$ & $\lambda_{N-1,0}$ & $\lambda_{N-1,0}E_{0,-1}%
+\lambda_{N-2,0}E_{0,-2}$\\\hline
$\vphantom{\vdots}1$ & $\lambda_{0,N-1}=\overline{\lambda_{N-1,0}}$ &
$\lambda_{0,N-1}E_{-1,0}+\lambda_{0,N-2}E_{-2,0}$\\\hline
\end{tabular}
\end{table}

Note that in non-generic cases, like $\lambda_{0}=1$, and/or ($\lambda
_{N-1}=1$ and $\lambda_{N-2}=0$), some of the listed eigenvectors are zero.
This is not surprising since these are exactly the cases where the root $1$ in
the characteristic polynomial is multiple. Let us consider these cases
separately:\setcounter{case}{1}

\begin{case}
\label{CasEigMult(1)}$\lambda_{0}=1$ ($\Rightarrow\lambda_{0,N-1}%
=\lambda_{0,N-2}=0$) and $\lambda_{N-1}<1$ or $\lambda_{N-2}>0$: In this case
the last two eigenvectors in the list in Table \textup{\ref{TabEigMult}} are
zero, but one may verify that these may be replaced by the vectors $E_{0,-2}$
and $E_{-2,0}$. Thus the eigenspace corresponding to $1$ is indeed
$2$-dimensional in this case, and spanned by $\openone_{\mathcal{K}}$ and
$E_{0,0}$.
\end{case}

\begin{case}
\label{CasEigMult(2)}$\lambda_{0}<1$, $\lambda_{N-1}=1$, $\lambda_{N-2}=0$: In
this case it follows from \textup{(\ref{eqRep.53})} that $\lambda_{N-1,k}=0$
for all $k\neq N-1$ and $\lambda_{N-2,k}=0$ for all $k\neq N-2$, so the matrix
in \textup{(\ref{eqRep.56})} degenerates into the upper-left-hand $3\times3$
matrix, which is%
\begin{equation}
\left(
\begin{array}
[c]{ccc}%
\lambda_{0} & 1-\lambda_{0} & 0\\
0 & 1 & 0\\
0 & 0 & 1
\end{array}
\right)  . \label{eqRep.57bis}%
\end{equation}
So we see that a new eigenvector for $\lambda_{N-1}-\lambda_{N-2}=1$ is%
\begin{equation}
E_{-2,-2}. \label{eqRep.57ter}%
\end{equation}
\end{case}

Let us now look at the even more degenerate case where the multiplicity of
$1$ in the characteristic polynomial is $3$:\setcounter{case}{0}

\begin{case}
\label{CasEigMult(3)}$\lambda_{0}=1$, $\lambda_{N-1}=1$, $\lambda_{N-2}=0$:
{}From the discussion of Case \textup{\ref{CasEigMult(2)}} above, it follows
that then the upper-left $3\times3$ matrix of \textup{(\ref{eqRep.56})} is
$\openone_{3}$, while the rest of the matrix is zero. Thus the fixed-point set
is the linear span of%
\begin{equation}
E_{0,0},\;E_{-1,-1}\text{, and }E_{-2,-2}. \label{eqRep.57tetra}%
\end{equation}
\end{case}

Thus the dimension of the eigensubspace is in all cases where the
parameters are restricted as in Remark \ref{RemRepNew.8} equal to the
multiplicity of the root $1$ in the characteristic polynomial, and the above
analysis implies

\begin{corollary}
\label{CorRep.7}If $g=2$ and $N>2$, the representation $T^{\left(  A\right)
}$ decomposes into at most three irreducible representations which are
mutually nonequivalent. Generically, $T^{\left(  A\right)  }$ is irreducible,
and otherwise the decomposition structure is summarized in Table
\textup{\ref{TabRepDecomp}.}
\end{corollary}

\begin{table}[ptbh]
\caption{Decompositions from eigenvalue-one multiplicities, and eigenvectors
for $\sigma^{\left(  A\right)  }$}%
\label{TabRepDecomp}
\addtolength{\tabcolsep}{-0.5\tabcolsep}
\begin{tabular}
[c]{l|c|l}%
\begin{tabular}
[t]{l}%
Case
\end{tabular}
&
\begin{tabular}
[t]{l}%
Number of\\
subrepresen-\\
tations
\end{tabular}
&
\begin{tabular}
[t]{l}%
The fixed-point set of $\sigma^{\left(  A\right)  }$\\
is spanned by
\end{tabular}
\\\hline%
\begin{tabular}
[t]{r}%
$\lambda_{0}=1$, $\lambda_{N-2}=0$,\\
\quad and $\lambda_{N-1}=1$%
\end{tabular}
& $3$ &
\begin{tabular}
[t]{l}%
$E_{0,0}$, $E_{-1,-1}$, and $E_{-2,-2}$%
\end{tabular}
\\\hline%
\begin{tabular}
[t]{r}%
$\lambda_{0}<1$, $\lambda_{N-1}=1$,\\
and $\lambda_{N-2}=0$%
\end{tabular}
& $2$ &
\begin{tabular}
[t]{l}%
$\openone_{K}$ and $E_{-2,-2}$\\
(i.e., $E_{0,0}+E_{-1,-1}$ and $E_{-2,-2}$)
\end{tabular}
\\\hline%
\begin{tabular}
[t]{l}%
$\lambda_{0}=1$ and\\
($\lambda_{N-1}<1$ or $\lambda_{N-2}>0$)
\end{tabular}
& $2$ &
\begin{tabular}
[t]{l}%
$\openone_{K}$ and $E_{0,0}$\\
(i.e., $E_{0,0}$ and $E_{-1,-1}+E_{-2,-2}$)
\end{tabular}
\\\hline
\end{tabular}
\end{table}

\begin{remark}
\label{RemRepNew.9}Table \textup{\ref{TabRepDecomp}} shows a remarkable
feature which these representations share with the $g=2$, $N=2$
representations studied in \cite{BEJ00}, namely that the fixed-point set of
$\sigma^{\left(  A\right)  }$ is both abelian and an algebra in all cases.
Hence it follows from \textup{(\ref{eqRep.12})} that the projection $P$ onto
$\mathcal{K}$ commutes with the commutant $\pi\left(  \mathcal{O}_{N}\right)
^{\prime}$, and thus we may find cyclic vectors for the various
subrepresentations by picking vectors in the range of the minimal projections
in the fixed-point set $\mathcal{B}\left(  \mathcal{K}\right)  ^{\sigma
^{\left(  A\right)  }}$. For example, for the first case in Table
\textup{\ref{TabRepDecomp}} \textup{(}Case \textup{\ref{CasEigMult(3)}),} the
three vectors $z^{-n}$, $n=0,1,2$, in $L^{2}\left(  \mathbb{T}\right)  $ are
cyclic for the three disjoint representations the original representation
decomposes into, respectively. For the second case in Table
\textup{\ref{TabRepDecomp}} \textup{(}Case \textup{\ref{CasEigMult(2)}),} the
pair $z^{0}$, $z^{-2}$ has the same property, as well as the pair $z^{-1}$,
$z^{-2}$ or any pair of the form $\lambda z^{0}+\mu z^{-1}$, $z^{-2}$. So far
we do not know a single example \textup{(}for general $N$, $g$\textup{)} where
$\mathcal{B}\left(  \mathcal{K}\right)  ^{\sigma^{\left(  A\right)  }}$ is not
abelian, or is not an algebra.
\end{remark}

\section{\label{Red}Reduction of representations as a reduction of $N$}

Consider an arbitrary element $A\in\mathcal{P}(\mathbb{T},\mathrm{U}(N))$.
Hence both $N$ and the genus $g$ are arbitrary, but given. We saw that there
is an associated wavelet filter $m^{(A)}$, as well as a representation
$T^{(A)}$ of $\mathcal{O}_{N}$ acting on $L^{2}(\mathbb{T})$. (As before we
work in the standard Fourier basis for $L^{2}(\mathbb{T})$, i.e.,
$e_{n}(z)=z^{n}$, $n\in\mathbb{Z}$.) We identified the subspace $\mathcal{K}$
spanned by $\{e_{-r},\ldots,e_{-1},e_{0}\}$ where $r$ was picked so that
reducibility of $T^{(A)}$ is decided by the fixed-point set of an associated
completely positive map $\sigma_{\mathcal{K}}^{(A)}$ given by
\begin{equation}
\sigma_{\mathcal{K}}^{(A)}(\,\cdot\,)=\sum_{i=0}^{N-1}V_{i}^{(A)}%
(\,\cdot\,)V_{i}^{(A)\ast} \label{eq6(1)}%
\end{equation}
with
\begin{equation}
V_{i}^{(A)}=\mathcal{P}_{\mathcal{K}}T_{i}^{(A)}\text{.} \label{eq6(2)}%
\end{equation}
While we identified the complete spectral picture of $\sigma_{\mathcal{K}%
}^{(A)}$ in the special case $g=2$, the higher-genus case is not yet entirely
understood. But we note that the identity from the $g=2$ case,
\begin{equation}
\sigma_{\mathcal{K}}^{(A)}(E_{0,0})=\lambda_{0}(A)E_{0,0}, \label{eq6(3)}%
\end{equation}
which is clear from (\ref{eqRep.56}), remains true for arbitrary $g$. Hence
our use of $\sigma_{\mathcal{K}}^{(A)}$ yields reducibility of the
representation $T^{(A)}$ whenever $\lambda_{0}(A)=1$. We recall that if
\begin{equation}
A(z)=A^{(0)}+zA^{(1)}+\cdots+z^{g-1}A^{(g-1)}, \label{eq6(4)}%
\end{equation}
then $\lambda_{0}(A)$ is defined as the $(0,0)$-matrix entry in $A^{(0)\ast
}A^{(0)}$.

\begin{remark}
\label{RemRedNew.1}To prove \textup{(\ref{eq6(3)})} for the general case,
recall that%
\begin{equation}
T_{i}^{\left(  A\right)  }e_{0}=\sum_{j=0}^{N-1}A_{i,j}\left(  z^{N}\right)
e_{j}=\sum_{j=0}^{N-1}\sum_{k=0}^{g-1}A_{i,j}^{\left(  k\right)  }e_{j+kN}.
\label{eq6(4)bis}%
\end{equation}
Using \textup{(\ref{eq6(2)}),} we therefore get%
\begin{equation}
V_{i}^{\left(  A\right)  }e_{0}=P_{\mathcal{K}}T_{i}^{\left(  A\right)  }%
e_{0}=\sum_{j=0}^{N-1}\sum_{k=0}^{g-1}A_{i,j}^{\left(  k\right)
}P_{\mathcal{K}}e_{j+kN}=A_{i,0}^{\left(  0\right)  }e_{0},
\label{eqRed.pound}%
\end{equation}
and by \textup{(\ref{eqRep.43}),}
\begin{align}
\sigma_{\mathcal{K}}^{\left(  A\right)  }\left(  E_{0,0}\right)   &
=\left( \sum_{i=0}^{N-1}\left|  A_{i,0}^{\left(  0\right)  }e_{0}\right\rangle
\left\langle A_{i,0}^{\left(  0\right)  }e_{0}\right| \right)  =\left(  \sum
_{i=0}^{N-1}\left|  A_{i,0}^{\left(  0\right)  }\right|  ^{2}\right)
\left|  e_{0}\right\rangle \left\langle e_{0}\right| \label{eqRed.poundbis}\\
&  =\left(  A^{\left(  0\right)  \,\ast}A^{\left(  0\right)  }\right)
_{0,0}E_{0,0}=\lambda_{0}\left(  A\right)  E_{0,0},\nonumber
\end{align}
which is \textup{(\ref{eq6(3)}).}

But it is worth stressing that the argument \textup{(\ref{eqRed.pound}),} for
the vector $e_{0}$ in $\mathcal{K}$, does \emph{not} carry over to the other
basis vectors $e_{-1}$, $e_{-2}$, $\dots$. This is clear already in the case
$g=2$ from \textup{(\ref{eqRep.56}).} More generally, for $e_{-1}$, for
example,
\begin{equation}
V_{i}^{\left(  A\right)  }e_{-1}=A_{i,0}^{\left(  1\right)  }e_{0}%
+A_{i,N-1}^{\left(  0\right)  }e_{-1}+A_{i,N-2}^{\left(  0\right)  }%
e_{-2}+\cdots, \label{eqRed.poundter}%
\end{equation}
and therefore%
\begin{multline}
\sigma_{\mathcal{K}}^{\left(  A\right)  }\left(  E_{-1,-1}\right)
=\lambda_{0,0}^{\left(  11\right)  }E_{0,0}+\lambda_{N-1,N-1}^{\left(
00\right)  }E_{-1,-1}+\lambda_{N-2,N-2}^{\left(  00\right)  }E_{-2,-2}%
+\cdots\label{eqRed.poundpound}\\
+\lambda_{N-1,0}^{\left(  01\right)  }E_{0,-1}+\lambda_{0,N-1}^{\left(
10\right)  }E_{-1,0}+\lambda_{N-2,N-1}^{\left(  00\right)  }E_{-1,-2}%
+\lambda_{N-1,N-2}^{\left(  00\right)  }E_{-2,-1}+\cdots,
\end{multline}
where $\lambda_{i,j}^{\left(  kl\right)  }:=\left(  A^{\left(  k\right)
\,\ast}A^{\left(  l\right)  }\right)  _{i,j}$, i.e., the $\left(  i,j\right)
$ entry in the $N\times N$ matrix $A^{\left(  k\right)  \,\ast}A^{\left(
l\right)  }$. Even if $\lambda_{N-1,N-1}^{\left(  00\right)  }=1$, that does
not imply that all the other coefficients $\lambda_{i,j}^{\left(  kl\right)
}$ from \textup{(\ref{eqRed.poundpound})} will necessarily vanish; see the
details below. While%
\begin{equation}
\lambda_{N-1,N-1}^{\left(  00\right)  }=1\Longrightarrow\lambda_{N-1,j}%
^{\left(  00\right)  }=\lambda_{i,N-1}^{\left(  00\right)  }=0\text{\qquad for
all }i,j\in\left\{  0,\dots,N-2\right\}  , \label{eqRed.disp2}%
\end{equation}
the other coefficients in \textup{(\ref{eqRed.poundpound})} such as
$\lambda_{N-2,N-2}^{\left(  00\right)  }$ or $\lambda_{0,0}^{\left(
11\right)  }$ would typically be nonzero even if $\lambda_{N-1,N-1}^{\left(
00\right)  }=1$. So even then, $E_{-1,-1}$ will not be fixed by $\sigma
_{\mathcal{K}}^{\left(  A\right)  }$, and there is not a natural analogue to
\textup{(\ref{eq6(3)}).} While if $g=2$, then the difference $\lambda
_{N-1,N-1}^{\left(  00\right)  }-\lambda_{N-2,N-2}^{\left(  00\right)  }$ is
in the spectrum of $\sigma_{\mathcal{K}}^{\left(  A\right)  }$, see Table
\textup{\ref{TabEigMult},} this will generally not be true if $g>2$. Then the
spectral picture for $\sigma_{\mathcal{K}}^{\left(  A\right)  }$ is not fully understood.
\end{remark}

However, the following result shows that a general wavelet representation
$T^{(A)}$ may be reducible because a special point $\lambda_{0}(A)$ in the
spectrum of $\sigma_{\mathcal{K}}^{(A)}$ may be one. This is a special
reduction, however, as we also showed that whenever any point $\lambda$ in
spectrum $\left(  \sigma_{\mathcal{K}}^{(A)}\right)  $ satisfies $\lambda=1$,
we get a possibly different reduction of the representation $(T^{(A)}%
,L^{2}(\mathbb{T}))$.

\begin{theorem}
\label{theorem6(1)}Let $A\in\mathcal{P}(\mathbb{T},\mathrm{U}(N))$. Then the
following five conditions are equivalent.

\begin{enumerate}
\item \label{theorem6(1)(1)}$\lambda_{0}(A)=1$.

\item \label{theorem6(1)(2)}There is a $B=(B_{i,j})_{i,j=1}^{N-1}%
\in\mathcal{P}(\mathbb{T},\mathrm{U}(N-1))$ such that
\begin{equation}
A(1)^{-1}A(z)=\left(
\begin{tabular}
[c]{c|cclc}%
$1$ & $0$ & $0$ & $\cdots$ & $0$\\\hline
$0$ & $B_{1,1}(z)$ & $\cdots$ & $\cdots$ & $B_{1,N-1}(z)$\\
$0$ & $\vdots$ &  &  & $\vdots$\\
$\vdots$ & $\vdots$ &  &  & $\vdots$\\
$0$ & $B_{N-1,1}(z)$ & $\cdots$ & $\cdots$ & $B_{N-1,N-1}(z)$%
\end{tabular}
\right)  , \label{eq6(5)}%
\end{equation}
i.e.,
\begin{equation}
A(1)^{-1}A(z)=\left(
\begin{tabular}
[c]{c|lll}%
$1$ & $0$ & $\cdots$ & $0$\\\hline
$0$ &  &  & \\
$\vdots$ &  & $B(z)$ & \\
$0$ &  &  &
\end{tabular}
\right)  =\left(  1\right)  \oplus B\left(  z\right)  . \label{eq6(5)bis}%
\end{equation}

\item \label{theorem6(1)(3)}After modifying with a $\mathrm{U}(N)$%
-automorphism of $\mathcal{O}_{N}$, the corresponding\linebreak wavelet filter
$m^{(A)}$ satisfies%
\begin{equation}
\left\{
\begin{array}
[c]{l}%
m_{0}^{(A)}(z)\equiv1\text{, for all }z\in\mathbb{T},\\
\vphantom{\vdots}m_{i}^{(A)}(z)=\sum_{j=1}^{N-1}B_{i,j}(z^{N})z^{j}.
\end{array}
\right.  \label{eq6(6)}%
\end{equation}

\item \label{theorem6(1)(4)}%
\[
T_{i}^{(A)^{\ast}}e_{0}\in\mathbb{C}e_{0}\text{, for all }i=0,\ldots
,N-1\text{.}%
\]

\item \label{theorem6(1)(5)}%
\[
\sigma_{\mathcal{K}}^{\left(  A\right)  \,\ast}\left(  E_{0,0}\right)
\in\mathbb{C}E_{0,0}.
\]
\end{enumerate}
\end{theorem}

\begin{remark}
\label{RemRedNew.3}Formula \textup{(\ref{eq6(6)})}
explains the assertion made in the Introduction about
reduction in the size of the number $N$ which serves as the scaling of the
wavelet in question. Once $B$ is identified then there is an $(N-1)$-wavelet
filter
\begin{equation}
m_{i-1}^{(B)}(z)=\sum_{j=0}^{N-2}B_{i,j+1}\left(  z^{N-1}\right)  z^{j},
\label{eq6(7)}%
\end{equation}
$i=1,\ldots,N-1$.
\end{remark}

\begin{proof}
\underline{(\ref{theorem6(1)(1)}) $\Rightarrow$ (\ref{theorem6(1)(2)})}: Using
the factorization (\ref{eqRep.22}) in Corollary \ref{CorRep.2}, we note that
\begin{multline}
A^{(0)\,\ast}A^{(0)}=\left(  \openone-Q_{g-1}\right)  \left(  \openone
-Q_{g-2}\right)  \cdots\label{eq6(7)bis}\\
\cdots\left(  \openone-Q_{1}\right)  \framebox{$\displaystyle\left
( \openone-Q_{0}\right) $}\left(  \openone-Q_{1}\right)  \cdots\\
\cdots\left(  \openone-Q_{g-2}\right)  \left(  \openone-Q_{g-1}\right)
\end{multline}
and the $V$ ($\in\mathrm{U}(N)$) from (\ref{eqRep.22}) is $V=A(1)$, i.e.,
evaluation of $A(z)$ at $z=1$. Let $\varepsilon_{0}$ be the first canonical
basis vector in $\mathbb{C}^{N}$. Then (\ref{theorem6(1)(1)}) states that
\begin{equation}
\left\|  \left(  \openone-Q_{0}\right)  \left(  \openone-Q_{1}\right)
\cdots\left(  \openone-Q_{g-1}\right)  \varepsilon_{0}\right\|  =1.
\label{eq6(7)ter}%
\end{equation}
Hence $\left\langle \left(  \openone-Q_{g-1}\right)  \varepsilon_{0}\mid
H\left(  \openone-Q_{g-1}\right)  \varepsilon_{0}\right\rangle =1$, where
\begin{equation}
H:=\left(  \openone-Q_{g-2}\right)  \cdots\left(  \openone-Q_{1}\right)
\framebox{$\displaystyle\left( \openone-Q_{0}\right) $}\left(  \openone
-Q_{1}\right)  \cdots\left(  \openone-Q_{g-2}\right)  . \label{eq6(7)tetra}%
\end{equation}
Using Schwarz's inequality and induction, we conclude that $\left(
\openone_{N}-Q_{j}\right)  \varepsilon_{0}=\varepsilon_{0}$, and therefore
$Q_{j}\varepsilon_{0}=0$ for all $j=0,\ldots,g-1$. Hence $\left(  \openone
_{N}-Q_{j}+zQ_{j}\right)  \varepsilon_{0}=\varepsilon_{0}$, and from
(\ref{eqRep.22}),
\begin{equation}
A(1)^{-1}A(z)\varepsilon_{0}=\varepsilon_{0}\text{.} \label{eq6(7)pent}%
\end{equation}
The same argument yields
\begin{equation}
\ip{\varepsilon_{0}}{A(1)^{-1}A(z)e_{j}}=\delta_{0,j}, \label{eq6(7)hex}%
\end{equation}
and so (\ref{theorem6(1)(2)}) follows.

\underline{(\ref{theorem6(1)(2)}) $\Rightarrow$ (\ref{theorem6(1)(3)})}:
Substituting (\ref{eq6(5)}) into (\ref{eq1.24}), we get
\begin{equation}
\left(  V^{-1}A\right)  _{i,0}(z)=\delta_{i,0}, \label{eq6(7)hept}%
\end{equation}
\begin{equation}
\left(  V^{-1}A\right)  _{0,j}(z)=\delta_{j,0}, \label{eq6(7)okt}%
\end{equation}
and
\begin{equation}
\left(  V^{-1}A\right)  _{i,j}(z)=B_{i,j}(z)\text{\qquad if }i,j\geq1,
\label{eq6(7)enn}%
\end{equation}
and (\ref{eq6(6)}) in (\ref{theorem6(1)(3)}) follows.

\underline{(\ref{theorem6(1)(3)}) $\Rightarrow$ (\ref{theorem6(1)(4)})}: Using%
\begin{equation}
T_{i}^{(A)\,\ast}\xi(z)=\frac{1}{N}\sum_{w^{N}=z}\overline{m_{i}^{(A)}(w)}%
\xi(w)\text{,\qquad}\xi\in L^{2}(\mathbb{T}), \label{eq6(8)}%
\end{equation}
and introducing the matrix product from (\ref{theorem6(1)(3)}), we get
\begin{equation}
\left\{
\begin{array}
[c]{l}%
A_{i,0}(z)=V_{i,0},\\
A_{i,j}(z)=\sum_{k=1}^{n-1}V_{i,k}B_{i,k}(z)\text{, }j>0.
\end{array}
\right.  \label{eq6(8)bis}%
\end{equation}
Using again (\ref{eq1.24}), we thus get
\begin{equation}
T_{i}^{(A)\,\ast}  \,e_{0}  =\overline{V_{i,0}}\,e_{0},
\label{eq6(8)ter}%
\end{equation}
where $V_{i,j}$ denotes the matrix entries of $V$ ($:=A(1)$) in $\mathrm{U}%
(N)$.

\underline{(\ref{theorem6(1)(4)}) $\Leftrightarrow$ (\ref{theorem6(1)(5)})
$\Rightarrow$ (\ref{theorem6(1)(1)})}: By (\ref{eq6(8)}),
\begin{equation}
T_{i}^{(A)\,\ast}e_{0}=\overline{A_{i,0}(z)}\in\mathcal{P}\left(
\mathbb{T},\mathbb{C}\right)  \text{;} \label{eq6(8)tetra}%
\end{equation}
and so if (\ref{theorem6(1)(4)}) is assumed, then
\begin{equation}
A_{i,0}^{(k)}=0\text{\qquad for }k>0. \label{eq6(8)pent}%
\end{equation}
Since
\begin{equation}
\sum_{i}\overline{A_{i,0}(z)}A_{i,0}(z)=1 \label{eq6(8)hex}%
\end{equation}
in general, we get
\begin{equation}
\lambda_{0}(A)=\sum_{i}\overline{A_{i,0}^{(0)}}A_{i,0}^{(0)}=1.
\label{eq6(8)hept}%
\end{equation}

{}From (\ref{eq6(8)tetra}), we get%
\begin{equation}
T_{i}^{\left(  A\right)  \,\ast}e_{0}=\sum_{k=0}^{g-1}\overline{A_{i,0}%
^{\left(  k\right)  }}\,e_{-k}, \label{eq6(8)okt}%
\end{equation}
and by (\ref{eqRep.44}),%
\begin{equation}
\sigma_{\mathcal{K}}^{\left(  A\right)  \,\ast}\left(  E_{0,0}\right)
=\sum_{k=0}^{g-1}\sum_{l=0}^{g-1}\lambda_{0,0}^{\left(  kl\right)  }E_{-k,-l}.
\label{theorem6(1)proof.pound}%
\end{equation}
Note the coefficient of $E_{0,0}$ in (\ref{theorem6(1)proof.pound}) is
$\lambda_{0}\left(  A\right)  =\lambda_{0,0}^{\left(  00\right)  }$. For the
other coefficients, we have%
\begin{equation}
\left|  \lambda_{0,0}^{\left(  kl\right)  }\right|  ^{2}\leq\lambda
_{0,0}^{\left(  kk\right)  }\lambda_{0,0}^{\left(  ll\right)  }
\label{theorem6(1)proof.poundbis}%
\end{equation}
and%
\begin{equation}
\lambda_{0,0}^{\left(  kk\right)  }=\sum_{i=0}^{N-1}\left|  A_{i,0}^{\left(
k\right)  }\right|  ^{2}. \label{theorem6(1)proof.poundter}%
\end{equation}
Hence, the vanishing of the non-$\left(  0,0\right)  $ coefficients in
(\ref{theorem6(1)proof.pound}) is equivalent to condition (\ref{eq6(8)pent}),
which we already showed is equivalent to both (\ref{theorem6(1)(4)}) and
(\ref{theorem6(1)(1)}).
\end{proof}

For $\xi \in \mathcal{K}$, let
$\mathcal{H}_{+}^{\left( A\right) }\left\lbrack \xi \right\rbrack $
denote the cyclic subspace in
$L^{2}\left( \mathbb{T}\right) $
generated by $\xi$ and the operators
$T_{i}^{\left( A\right) }$,
and let 
$E_{+}^{\left( A\right) }\left\lbrack \xi \right\rbrack $
be the corresponding projection.

\begin{corollary}
\label{CorRedNew.4}Suppose condition
\textup{(\ref{theorem6(1)(1)})} of Theorem \textup{\ref{theorem6(1)}} holds.
Then the \textup{(}unique\/\textup{)} operator $W$
in the commutant of 
$T^{\left( A\right) }$
which satisfies
\begin{equation}
P_{\mathcal{K}}WP_{\mathcal{K}}=E_{0,0} \label{eqCorRedNew.4pound}
\end{equation}
is $W=E_{+}^{\left( A\right) }\left\lbrack e_{0} \right\rbrack $.
\end{corollary}

\begin{proof}
Let $v\in\mathbb{C}^{N}$ be the
vector with coordinates
$v_{i}:=A_{i,0}^{\left( 0\right) }$. Then (\ref{theorem6(1)(1)}) holds
if and only if $\left\| v\right\| =1$. In
that case, the restriction of 
$T^{\left( A\right) }$
to $\mathcal{H}_{+}^{\left( A\right) }\left\lbrack e_{0} \right\rbrack $
is the subrepresentation
which is induced by the Cuntz
state \cite{Cun77} $\omega_{v}$. Specifically,
if $I=\left( i_{1},\dots ,i_{k}\right) $,
$J=\left( j_{1},\dots ,j_{l}\right) $
are multi-indices, then it follows
from (\ref{theorem6(1)(4)}) that
\begin{equation}
\ip{e_{0}}{T_{I}^{\left( A\right) }T_{J}^{\left( A\right) \,\ast}e_{0}}
=v_{i_{1}}\cdots v_{i_{k}}\overline{v_{j_{1}}}\cdots \overline{v_{j_{l}}}
=\omega_{v}\left( s_{I}^{}s_{J}^{\ast}\right) .
\label{eqCorRedNew.4proofpoundpound}
\end{equation}
Hence 
$\mathcal{H}_{+}^{\left( A\right) }\left\lbrack e_{0} \right\rbrack $
is the closed subspace spanned by $e_{0}$ and
$\left\{ T_{I}^{\left( A\right) }e_{0}\right\} $.
Since
\begin{equation}
T_{i}^{\left( A\right) }e_{0}\left( z\right) 
=\sum_{j=0}^{N-1}A_{i,j}\left( z^{N}\right) z^{j}
\in\mathcal{P}\left( \mathbb{T},\mathbb{C}\right) ,
\label{eqCorRedNew.4proofpoundpoundbis}
\end{equation}
$\mathcal{H}_{+}^{\left( A\right) }\left\lbrack e_{0} \right\rbrack $
is contained
in the Hardy space
$\mathcal{H}_{+}
:=\overline{\operatorname*{span}}\left\{ e_{n};n\geq 0\right\} $.
Since
$\mathcal{K}
=\operatorname*{span}\left\{ e_{-r},\dots ,e_{-1},e_{0}\right\} $,
it follows that
$E_{+}^{\left( A\right) }\left\lbrack e_{0} \right\rbrack $
satisfies (\ref{eqCorRedNew.4pound}). Since
$\sigma^{\left( A\right) }\left( E_{0,0}\right) =E_{0,0}$,
the operator $W$ in
the commutant of $T^{\left( A\right) }$ satisfying
(\ref{eqCorRedNew.4pound}) is unique by Theorem \ref{ThmRep.1},
and so
$W=E_{+}^{\left( A\right) }\left\lbrack e_{0} \right\rbrack $.
\end{proof}

\begin{remark}
\label{RemRedNew.5}The
Hilbert space prior to
Corollary \textup{\ref{CorRedNew.4}}
is of the form
$\mathcal{H}_{+}^{\left( A\right) }\left\lbrack \xi \right\rbrack $
for some $\xi\in\mathcal{K}$ where
\begin{equation}
T_{i}^{\left( A\right) \,\ast }\xi =v_{i}\xi 
\label{eqCorRedNew.4proofpoundpoundpound}
\end{equation}
for some 
$v=\left( v_{i}\right) \in \mathbb{C}^{N}$.
Since
$\omega_{v}$ is a Cuntz state, the
corresponding
representation is unique from $v$ up
to unitary equivalence, by
\cite{Cun77}. It also follows
from \cite[Theorem 6.3]{Jor99a} that
every $\xi \in \mathcal{K}$ which
satisfies \textup{(\ref{eqCorRedNew.4proofpoundpoundpound})} for some
$v\in \mathbb{C}^{N}$ must be a monomial,
i.e., $\xi \left( z\right) =z^{-k}$ for some
$k\in \left\{ 0,1,\dots ,r\right\} $.
In our special setting, this can
be proved as follows: If $T_{i}^{\ast}\xi =\bar{v}_{i}\xi $
for a unit vector $\xi \in \mathcal{K}$, where
$\sum_{i}\left| v_{i}\right| ^{2}=1$, then
$\xi =\sum_{i}T_{i}^{}T_{i}^{\ast}\xi 
=\sum_{i}\bar{v}_{i}T_{i}\xi $, or, spelled out,
$\xi \left( z\right) 
=\sum_{i}\bar{v}_{i}m_{i}\left( z\right) \xi \left( z^{N}\right) $. Hence,
putting $m\left( z\right) =\sum_{i}\bar{v}_{i}m_{i}\left( z\right) $, we have
$\xi \left( z\right) =m\left( z\right) \xi \left( z^{N}\right) $. Now we may
use the argument from \cite[page 104]{Jor99a}
or \cite[Theorem 3.1]{BrJo97b} to conclude
that $\left| m\left( z\right) \right| =1$ and
$\left| \xi \left( z\right) \right| =1$ for almost
all $z\in \mathbb{T}$. But both $m\left( z\right) $ and
$\overline{\xi \left( z\right) }$
are polynomials, and we deduce from
Lemma \textup{\ref{lemma2.1}} that they are monomials.
If $\xi \left( z\right) =z^{-k}$ it follows that
$m\left( z\right) =\xi \left( z\right) \overline{\xi \left( z^{N}\right) }
=z^{\left( N-1\right) k}$. Thus
a very restrictive necessary and sufficient
condition for the Cuntz state $\omega _{v}$
to occur is that
\begin{equation}
\sum_{i}\bar{v}_{i}m_{i}\left( z\right) 
=z^{\left( N-1\right) k}
\text{\qquad for some }
k\in \left\{ 0,1,\dots ,r\right\} . \label{eqRemRedNew.5(1)}
\end{equation}

The
other
details of the decompositions may be spelled out as follows:
Let $\mathcal{H}=L^{2}\left( \mathbb{T}\right) $.
Note first that Corollary
\textup{\ref{CorRedNew.4}} holds whenever $E_{0,0}$ is replaced
by any one-dimensional projection
in the fixed-point set of $\sigma$
in the generality of Theorem \textup{\ref{ThmRep.1},}
since
$E_{0,0}\left( \mathcal{K}\right) =E_{0,0}\left( \mathcal{H}\right) $
is then a
one-dimensional $T^{\left( A\right) \,\ast }$-invariant space.
If $e$ is a general
\textup{(}not necessarily one-dimensional\/\textup{)}
projection in the
fixed-point set, it follows from Lemma
\textup{3.3} in \cite{BJKW00} that $e$ commutes
with all the operators
$PT_{i}^{\left( A\right) \,\ast }P$ and $PT_{i}^{\left( A\right) }P$
and hence the projection $E$ onto
$\left\lbrack \mathcal{O}_{N}e\mathcal{K}\right\rbrack $
in $\mathcal{H}$ is in the commutant of the representation,
and $PEP=e$. This justifies claims
in Remark \textup{\ref{RemRepNew.9}.} In the case $g=2$, which
is
completely enumerated in
\cite{BEJ00}
\textup{(}for $N=2$\textup{)} and in Table \textup{\ref{TabRepDecomp}}
\textup{(}for $N>2$\textup{)}, the ranges of the
projections in the fixed-point set
are spanned by subsets of the
orthonormal set
$\left\{ e_{0},e_{-1},\dots ,e_{-r}\right\} $,
and hence each Hilbert space in the decomposition
has the form $E_{+}^{\left( A\right) }\left\lbrack \xi \right\rbrack $ for
a suitable $\xi $ in this set.
The vector $\xi $ defines a Cuntz
state if and only if the corresponding
projection in
$\mathcal{B}\left( \mathcal{K}\right) ^{\sigma}$ is
one-dimensional.
Note also that the result \cite[Theorem 6.3]{Jor99a}
spelled out above implies that all one-dimensional
projections in $\mathcal{B}\left( \mathcal{K}\right) ^{\sigma}$
are diagonal in the standard basis.
This gives a partial explanation of these results.
\end{remark}

\subsection{Other examples of decompositions of
$L^{2}\left( \mathbb{T}\right) $}

Generically, the wavelet examples
do not have
minimal projections $e$ in
$\mathcal{B}\left( \mathcal{K}\right) ^{\sigma}$
of dimension one,
and even if there is one such
projection in
$\mathcal{B}\left( \mathcal{K}\right) ^{\sigma}$,
there may be others that do not have dimension one.
For illustration, let us repeat
some of the examples in \cite{BEJ00},
rephrased in the setting of the present paper.
Take, for example,
$N=2=g$. Then, by \textup{(\ref{eqRemRedNew.5(1)}),} only the two cases
$
\begin{pmatrix}
1 & 0 \\
0 & z
\end{pmatrix}
$ and
$
\begin{pmatrix}
z & 0 \\
0 & 1
\end{pmatrix}
$ have the Cuntz-state vectors. For
$A\left( z\right) =
\begin{pmatrix}
1 & 0 \\
0 & z
\end{pmatrix}
$,
even though $\lambda_{0}\left( A\right) =1$,
this example has a
minimal $2$-dimensional
projection
$e$ in
$\mathcal{B}\left( \mathcal{K}\right) ^{\sigma}$,
namely $e=E_{-1,-1}+E_{-2,-2}$.
It is minimal in the sense that
the restriction of $T^{\left( A\right) }$
to $E_{e}=\left\lbrack \mathcal{O}_{2}e\mathcal{K}\right\rbrack $
is irreducible; see \cite[(4.52)]{BEJ00}.
The matrix $A\left( z\right) =
\begin{pmatrix}
0 & 1 \\
z & 0
\end{pmatrix}
$ has $\lambda_{0}\left( A\right) =0$, and its representation
$T^{\left( A\right) }$ decomposes into a sum of two
irreducibles
associated with respective $2$-dimensional projections $e$ and $f$
in $\mathcal{B}\left( \mathcal{K}\right) ^{\sigma}$:
$e=E_{0,0}+E_{-1,-1}$ and $f=E_{-2,-2}+E_{-3,-3}$; see \cite[(4.46)]{BEJ00}.
Hence in this case,
$\mathcal{H}=L^{2}\left( \mathbb{T}\right) 
=\left\lbrack \mathcal{O}_{2}e\mathcal{K}\right\rbrack 
\oplus \left\lbrack \mathcal{O}_{2}f\mathcal{K}\right\rbrack $
with $T^{\left( A\right) }$ restricting to an
irreducible representation on each
of the two subspaces.
A direct application of (\ref{eqRep.14})--(\ref{eqRep.15})
in Theorem \ref{ThmRep.1} shows that these
two subrepresentations are
disjoint, i.e., inequivalent.
While these
examples have $g=2$, the cases $g>2$
entail a richer decomposition structure.
\medskip

Note that the result, Theorem \textup{\ref{ThmRep.1},}
which is used in analyzing the
representations $T^{\left( A\right) }$, applies also when
$g>2$. Moreover, the factorization result,
Corollary \textup{\ref{CorRep.2},} yields in principle
a way of understanding the general
case, i.e., arbitrary $g$ and $N$, but unpublished
experimentation with examples for
$g=3$ \textup{(}i.e., two projections in $\mathbb{C}^{N}$\textup{)}
has not so far yielded decomposition
structures more general than
the above mentioned ones.

\begin{remark}
\label{RemRedNew.4}
The proof of \textup{(\ref{eq6(3)})} shows more
generally that
\begin{equation}
\sigma^{\left( B,A\right) }\left( E_{0,0}\right) 
=\left( \sum_{i=0}^{N-1}
\overline{A_{i,0}^{\left( 0\right) }}
B_{i,0}^{\left( 0\right) }\right) E_{0,0}.
\label{eqRedNew.34}
\end{equation}
Referring to Theorem \textup{\ref{ThmRep.1},} we note then that
the necessary and sufficient condition for $E_{0,0}$ to
induce an operator in $L^{2}\left( \mathbb{T}\right) $ which intertwines
the two representations $T^{\left( A\right) }$ and $T^{\left( B\right) }$ is
$\left( A^{\left( 0\right) \,\ast}B^{\left( 0\right) }\right) _{0,0}=1$.
Moreover, if
this holds, then both $T^{\left( A\right) }$ and $T^{\left( B\right) }$
must satisfy the equivalent conditions in
Theorem \textup{\ref{theorem6(1)}}, and we must then
further have $A_{i,0}^{\left( 0\right) }=B_{i,0}^{\left( 0\right) }$
for all $i$.
Then the intertwining
operator $W$ which
is induced by $E_{0,0}$
via 
$\sigma^{\left( B,A\right) }\left( E_{0,0}\right) =E_{0,0}$
is the one which
results from the uniqueness
of the $\omega_{v}$-representation
where $v_{i}=A_{i,0}^{\left( 0\right) }
=B_{i,0}^{\left( 0\right) }$.
\end{remark}

\begin{acknowledgements}
This research was done when one of us \textup{(}P.E.T.J.\textup{)} visited the
University of Oslo with support from the university and NFR. He is grateful
for kind hospitality. We wish to thank Brian Treadway for his skillful
typesetting, and Rune Kleveland and Brian Treadway for suggestions
and elimination of mistakes in earlier versions of this paper.
We are also indebted to Rong-Qing Jia and the referee
for bringing the references
\cite{JiMi91},
\cite{JiSh94},
\cite{LLS96},
\cite{RiSh91},
\cite{RiSh92}
to our attention, and thus enabling us to complete
the discussion between (\ref{eq3.11}) and (\ref{eq3.11bis}).

This is an expanded version of the invited lecture delivered by P.E.T.J. at
the
International Conference on Wavelet Analysis and Its Applications
held at
Zhongshan University, Guangzhou, China,
in November of 1999
(the conference announcement may be seen in
\textit{The Wavelet Digest} \textbf{8} (1999), on the World Wide Web at
{http://www.wavelet.org/wavelet/digest\_08/digest\_08.04.html\#13}).
\end{acknowledgements}

%\bibliographystyle{bftalpha}
%\bibliography{jorgen}

\end{document}